\documentclass[]{amsart}
\usepackage{hyperref}
\usepackage{CarnotHeader}
\theoremstyle{plain}

\usepackage{yfonts}

\newcommand{\HH}{\mathbb{H}}
\newcommand{\scr}[1]{\mathscr{#1}}
\newcommand{\frk}[1]{\mathfrak{#1}}
\newcommand{\cal}[1]{\mathcal{#1}}
\newcommand{\W}{\mathbb{W}}
\newcommand{\V}{\mathbb{V}}
\newcommand{\Span}{\mathrm{span}}

\renewcommand{\G}{\mathbb{G}}
\newcommand{\g}{\mathfrak{g}}

\newcommand{\ad}{\mathtt{ad}}
\newcommand{\Ad}{\mathtt{Ad}} 
\newcommand{\PP}{\mathbb{P}}
\newcommand{\MM}{\mathbb{M}}
\newcommand{\LL}{\mathbb{L}}
\newcommand{\KK}{\mathbb{K}}

\renewcommand{\DH}{D_H\;\!\!}

\newcommand{\aveint}[2]{\mathchoice%
          {\mathop{\kern 0.2em\vrule width 0.6em height 0.69678ex depth -0.58065ex
                  \kern -0.8em \intop}\nolimits_{\kern -0.45em#1}^{#2}}%
          {\mathop{\kern 0.1em\vrule width 0.5em height 0.69678ex depth -0.60387ex
                  \kern -0.6em \intop}\nolimits_{#1}^{#2}}%
          {\mathop{\kern 0.1em\vrule width 0.5em height 0.69678ex depth -0.60387ex
                  \kern -0.6em \intop}\nolimits_{#1}^{#2}}%
          {\mathop{\kern 0.1em\vrule width 0.5em height 0.69678ex depth -0.60387ex
                  \kern -0.6em \intop}\nolimits_{#1}^{#2}}}

\title[Area and coarea in Carnot Groups]{Area of intrinsic graphs and coarea formula\\ in Carnot groups}

\author[Julia]{Antoine Julia}
\author[Nicolussi Golo]{Sebastiano Nicolussi Golo}
\author[Vittone]{Davide Vittone}
\address[Julia and Nicolussi Golo and Vittone]{Dipartimento di Matematica ``T.Levi-Civita'', via Trieste 63, 35121 Padova, Italy.}
\email{antoine.julia@math.unipd.it}
\email{sebastiano.nicolussi@math.unipd.it}
\email{vittone@math.unipd.it}
\thanks{The authors are supported by the University of Padova STARS Project ``Sub-Riemannian Geometry and Geometric Measure Theory Issues: Old and New'' (SUGGESTION). They also acknowledge the support of GNAMPA of INdAM and   FFABR 2017 of MIUR (Italy).}

\date{\today}



\begin{document}
\begin{abstract}
We consider submanifolds of sub-Riemannian Carnot groups with intrinsic $C^1$ regularity ($C^1_H$). Our first main result is an area formula for $C^1_H$ intrinsic graphs; as an application, we deduce density properties for Hausdorff measures on rectifiable sets. Our second main result is a coarea formula for slicing $C^1_H$ submanifolds into level sets of a $C^1_H$ function.
\end{abstract}

\maketitle
\section{Introduction}
The interest
towards Analysis and Geometry in Metric Spaces  grew drastically in the last decades: a major effort has been devoted to the development of analytical tools for the study of geometric problems, and sub-Riemannian Geometry provided a particularly fruitful setting for these investigations. The present paper aims at giving a contribution in this direction by providing some geometric integration formulae, namely: an {\it area formula} for submanifolds with (intrinsic) $C^1$ regularity, and a {\it coarea formula} for slicing such  submanifolds into level sets of maps with (intrinsic) $C^1$ regularity.

We will work in the setting of a {\it Carnot group} $\G$, i.e., a connected, simply connected and nilpotent Lie group with stratified Lie algebra. We refer to  Section~\ref{sec12101536} for precise definitions; here, we only recall that Carnot groups have a distinguished role in sub-Riemannian Geometry, as they provide the infinitesimal models (tangents spaces) of sub-Riemannian manifolds, see e.g.~\cite{Bellaiche}. As usual, a Carnot group is endowed with a  distance $\rho$ that is left-invariant and 1-homogeneous with respect to the group dilations. 

Our main objects of investigation are {\it $C^1_H$ submanifolds}, which are introduced as (noncritical) level sets of functions with intrinsic $C^1$ regularity: let us briefly introduce the relevant definitions, which are  more precisely stated in Section~\ref{sec:preliminaries}.  Given an open set $\Omega\subset\G$ and another Carnot\footnote{One could more generally assume that $\G'$ is only graded, see Remark~\ref{rem:rem2.3}.} group $\G'$, a map $f:\Omega\to\G'$ is said to be {\it of class $C^1_H$} if it is  differentiable \`a la P.~Pansu~\cite{Pansu} at all $p\in\Omega$ and the differential $\DH f_p:\G\to\G'$ is continuous in $p$. Let us mention that the $C^1_H$ regularity of $f$ is equivalent to its {\it strict} Pansu differentiability (see Proposition~\ref{prop:C1Hiffstrict}): such a notion is introduced in Section~\ref{sec:Pansudifferentiability} and turns out to be useful for simplifying several arguments. Given a Carnot group $\G'$, a set $\Sigma\subset\G$ is a $C^1_H(\G;\G')$-submanifold if it is locally a level set of a map $f:\G\to\G'$  of class $C^1_H$  such that, at all points $p$, $\DH f_p$ is surjective and $\ker \DH  f_p$  splits $\G$. We say that a normal homogeneous subgroup $\W<\G$ {\it splits} $\G$ if there exists another homogeneous subgroup $\V<\G$, which is {\it complementary} to $\W$, i.e.,  such that $\V\cap \W = \{0\}$ and $\G= \W\V$. Observe that $\V$ is necessarily isomorphic to $\G'$, see Remark~\ref{rem:V=G'}. We will also say that $p$ is {\it split-regular} for $f$ if  $\DH f_p$ is surjective and $\ker \DH  f_p$ splits $\G$.

In Sections~\ref{sec09181614} and~\ref{sec:C1H_and_rectifiable} we prove that an Implicit Function Theorem holds for a $C^1_H$ submanifold $\Sigma$; namely, $\Sigma$ is (locally) an {\it intrinsic graph}, i.e.,  there exist complementary homogeneous subgroups $\W,\V$ of $\G$  and a function $\phi:A\to\V$ defined on an open subset $A\subset\W$ such that $\Sigma$ coincides with the intrinsic graph $\{w\phi(w):w\in A\}$ of $\phi$. The function $\phi$ is of class $C^1_{\W,\V}$ (see Definition~\ref{def:C1WV}) and it turns out to be {\it intrinsic Lipschitz continuous} according to the theory developed in recent years by B.~Franchi, R.~Serapioni and F.~Serra~Cassano, see e.g.~\cite{FrSeSe2006intrinsic,FSSCJGA2011,FranchiSerapioni2016Intrinsic}. We have to mention that both the Implicit Function Theorem and the intrinsic Lipschitz continuity of $\phi$ follow also from~\cite[Theorem 1.4]{MR3123745}: the proofs we provide in Sections~\ref{sec09181614}--\ref{sec:C1H_and_rectifiable}, however, seem shorter than those in~\cite{MR3123745} and allow for some finer results we need, see e.g.~Lemmas~\ref{lem09122015} and~\ref{lem09151054}. For related results, see~\cite{ArenaSerapioni,CittiManfredini,FSSCCAG,FSSCAdvMath,VSNS}.

Our first main result is an area formula for intrinsic graphs of class $C^1_{\W,\V}$ (hence, in particular, for $C^1_H$ submanifolds) where complementary subgroups $\W<\G$ and $\V<\G$ are fixed with $\W$ normal. Throughout the paper we denote by $\psi^d$ either the spherical or the Hausdorff measure of dimension $d$ in $\G$.

\begin{Theorem}[Area formula]\label{prop05161504}
	Let $\G$ be a Carnot group 
	and let $\G=\W\V$ be a splitting.
	Let $A\subset\W$ be an open set, $\phi\in C^1_{\W,\V}(A)$ and let $\Sigma:=\{w\phi(w):w\in A\}$ be the intrinsic graph of $\phi$; let $d$ be  the homogeneous dimension of $\W$. 
	Then, for all Borel functions $h:\Sigma\to[0,+\infty)$,
	\begin{equation}\label{eq09071257}
	\int_\Sigma h \dd \psi^{d} 
	= \int_A h(w\phi(w)) \calA(T_{w\phi(w)}^H\Sigma) \dd \psi^d (w) .
	\end{equation}
\end{Theorem}

The function $\calA(\,\cdot\,)$ appearing in~\eqref{eq09071257} is continuous and it is called {\it area factor}: it is defined in Lemma~\ref{lem:areafactor} and it depends  only on ($\W,\V$ and) the {\it homogeneous tangent space} $T^H_p\Sigma$ at points $p\in\Sigma$. The definition of area factor in Lemma~\ref{lem:areafactor} is only implicit, but of course we expect it can be made more explicit in terms of suitable derivatives of the map $\phi$: to the best of our knowledge, this program has been completed only in Heisenberg groups, see e.g.~\cite{AmbSerVit2006Intrinsic,ArenaSerapioni,Corni,CorniMagnani,FSSCAdvMath}. A relevant tool in the proof of Theorem~\ref{prop05161504} is a differentiation theorem for measures (Proposition~\ref{propFedDensity}) which is based on the so-called {\it Federer density}~\eqref{eq:defFedererdensity}: the importance of this notion was pointed out only recently by V.~Magnani, see~\cite{MagnaniEdinburgh,MagnaniNewDifferentiation,MagnaniSomeRemarks} and~\cite{FSSCcentered}. Observe that the validity of a (currently unavailable) Rademacher-type Theorem for intrinsic Lipschitz graphs would likely allow to extend Theorem~\ref{prop05161504} to the case of intrinsic Lipschitz $\phi$. 

A first interesting consequence of Theorem~\ref{prop05161504} is the following Corollary~\ref{cor:SvsH}, which is reminiscent of  the well-known equality between  Hausdorff and spherical Hausdorff measures on $C^1$ submanifolds (and, more generally, on rectifiable subsets) of $\R^n$. We refer to Definitions~\ref{def:rectifiable} and~\ref{def:apprtangent} for the notions of {\it countably $(\G;\G')$-rectifiable} set $R\subset\G$ and of {\it approximate tangent space} $T^H R$. 
Such sets have Hausdorff dimension $Q-m$, where $Q$ and $m$ denote, respectively, the homogeneous dimensions of $\G$, $\G'$; we write  $\cal H^{Q-m},\cal S^{Q-m}$, respectively, for Hausdorff and spherical Hausdorff measures. We denote by $\scr T_{\G,\G'}$ the space of possible tangent subgroups to $(\G;\G')$-rectifiable sets\footnote{Equivalently, $\scr T_{\G,\G'}$ is the space of normal subgroups $\PP<\G$ for which there exist a complementary subgroup in $\G$ and a surjective homogeneous morphism $L:\G\to\G'$ such that $\PP=\ker L$.} and, by abuse of notation, we  write $T^HR$ for the map $R\ni p\mapsto T^H_p R\in\scr T_{\G,\G'}$.

\begin{Corollary}\label{cor:SvsH}
Let $\G,\G'$ be Carnot groups of homogeneous dimensions $Q$, $m$, respectively. Then, there exists a continuous function $\textfrak a:\scr T_{\G,\G'}\to [1,2^{Q-m}]$ such that, for every countably $(\G;\G')$-rectifiable set $R\subset\G$
\begin{equation}\label{eq:AAAAA}
\cal S^{Q-m} \hel R= \textfrak a(T^H R)\cal H^{Q-m}\hel R\,.
\end{equation}
Moreover, if $\G$ is a Heisenberg group $\HH^n$ with a rotationally invariant distance $\rho$ and $\G'=\R$, then the function $\textfrak a$ is constant, i.e., there exists $C\in[1,2^{2n+1}]$ such that 
\[
\text{$\cal S^{2n+1} \hel R= C\cal H^{2n+1}\hel R\qquad \forall\:(\HH^n,\R)$-rectifiable set $R\subset\HH^n$.}
\]
\end{Corollary}

Heisenberg groups and rotationally invariant distances are defined in Section~\ref{sec12101536} by condition \eqref{eq:seba1}, while Corollary~\ref{cor:SvsH} is proved in Section~\ref{sec:area}. To the best of our knowledge, this result is new even in the first Heisenberg group $\HH^1$, see also~\cite[page 359]{MagnaniSomeRemarks}. Corollary~\ref{cor:SvsH} is deeply connected to the {\it isodiametric problem}, see Remark~\ref{rem:isodiametric}.

Not unrelated with Corollary~\ref{cor:SvsH} is another interesting consequence of Theorem~\ref{prop05161504}, namely, the existence of the {\it density} of Hausdorff and spherical measures on rectifiable sets. In Corollary~\ref{cor:densityexistence} we indeed prove that, if $R\subset\G$ is $(\G;\G')$-rectifiable, then the limit 
\[
\textfrak d(p):=\lim_{r\to 0^+}\frac{\psi^{Q-m}(R\cap \Ball(p,r))}{r^{Q-m}}
\]
exists for $\psi^{Q-m}$-a.e.~$p\in R$, where $\Ball(p,r)$ is the open ball of center $p$ and radius $r$ for the distance of $\G$. Actually, $\textfrak d(p)$ depends only on $T^H_pR$, in a continuous way. When $\G$ is the Heisenberg group $\HH^n$ endowed with a rotationally invariant distance, $\G'=\R^m$ for some $1\le m\le n$, and $\psi$ is the spherical measure, then $\textfrak d$ is constant, see Corollary~\ref{cor:densityconstantinHeis}.

The area formula is a key tool also in the proof of our second main result, the coarea formula in Theorem~\ref{thm:coarea} below. The classical coarea formula was first proved in the seminal paper~\cite{FedererCurvatureMeasures} and it is one of the milestones of Geometric Measure Theory. Sub-Riemannian coarea formulae have been obtained in \cite{MagnaniPublMat,MagnaniMathNacr,MagnaniCEJM,MagnaniNonHorizontal, zbMATH06358560, zbMATH06235931}, assuming classical (Euclidean) regularity on the slicing function $u$, and in \cite{MagnaniAreaCoarea,MagnaniStepanovTrevisan,MontiVittoneHeight},  assuming intrinsic regularity but only in the setting of the Heisenberg group. Here we try to work in the utmost generality: we consider a $C^1_H$ submanifold $\Sigma\subset\G$, seen as the level set of a $C^1_H$ map $f$ with values in a homogeneous group $\MM$, and we {\it slice} it into level sets of a map $u$ with values into another homogeneous group $\LL$. We assume for the sake of generality (see below) that $\LL,\MM$ are complementary subgroups of a larger homogeneous group $\KK=\LL\MM$; we also denote by $Q,m,\ell$ the homogeneous dimensions of $\G,\MM,\LL$, respectively.

\begin{Theorem}[Coarea formula]\label{thm:coarea}
	Let $\G,\LL,\MM$ be Carnot groups and let $\Omega\subset\G$ be open. Fix $f\in C^1_H(\Omega;\MM)$ and  assume that all points in $\Omega$ are split-regular for $f$, so that  $\Sigma:=\{p\in\Omega:f(p)=0\}$ is a $C^1_H$ submanifold. Consider a function $u:\Omega\to\LL$  such that $uf\in C^1_H(\Omega;\KK)$ and assume that
\begin{equation}\label{eq:technicalassumption}
\text{for $\psi^{Q-m}$-a.e.~$p\in\Sigma$,}\quad
\left\{
\begin{array}{l}
\text{either $\DH  (uf)_p|_{T^H_p\Sigma}$ is not surjective on $\LL$,}\\
\text{or $p$ is split-regular for $uf$.}
\end{array}
\right.
\end{equation}		
	Then, for every Borel function $h:\Sigma\to[0,+\infty)$ the  equality
	\begin{equation}\label{eq:coarea2}
	\int_\Sigma h(p)\,\calC (T^H_p\Sigma,\DH  (uf)_p) \, \dd\psi^{Q-m} (p)=\!\int_\LL \int_{\Sigma\cap u^{-1}(s)} h(p)\dd\psi^{Q-m-\ell}(p)\:\dd\psi^\ell(s) 
		\end{equation}
	holds.
\end{Theorem}

In~\eqref{eq:coarea2}, the symbol $\calC (T^H_p\Sigma,\DH  (uf)_p)$ denotes the {\it coarea factor}: let us stress that  it depends only on the restriction of $u$ to $\Sigma$ and that it does not depend on the choice of $f$ outside of $\Sigma$, see Remark~\ref{rem:olyonu}. The $\psi^\ell$-measurability of the function $\LL\ni s\mapsto\int_{\Sigma\cap u^{-1}(s)} h\,\dd\psi^{Q-m-\ell}$ is part of the statement.

The assumption $uf\in C^1_H(\Omega;\KK)$ becomes more transparent when $\KK=\LL\times\MM$ is a direct product (roughly speaking, when  $\LL,\MM$ are ``unrelated'' groups): in this case, it is in fact equivalent to the $C^1_H$ regularity of $u$. 
Moreover, since $T^H_p\Sigma=\ker \DH  f_p$,  the equality $\DH (uf)_p|_{T^H_p\Sigma}=\DH  u_p|_{T^H_p\Sigma}$ holds.
Eventually, the statement of Theorem~\ref{thm:coarea} can  at the same time be simplified, stated in a more natural way, and generalized to rectifiable sets, as follows.

\begin{Corollary}\label{cor:coareaPRODUCT}
Let $\G, \LL,\MM$ be Carnot groups, let $\Omega\subset\G$ be an open set and let $R\subset\Omega$ be $(\G;\MM)$-rectifiable; assume that  $u\in C^1_H(\Omega;\LL)$ is such that
\begin{equation}\label{eq:technicalassumptionPRODUCT}
\text{for $\psi^{Q-m}$-a.e.~$p\in R$,}\quad
\left\{
\begin{array}{l}
\text{either $\DH  u_p|_{T^H_pR}$ is not surjective on $\LL$,}\\
\text{or $T^H_pR\cap\ker\DH  u_p$ splits $\G$.}
\end{array}
\right.
\end{equation}
Then, for every Borel function $h:\Omega\to[0,+\infty)$ the  equality
	\begin{equation}\nonumber 
	\int_R h(p)\,\calC (T^H_p R,\DH  u_p) \, \dd\psi^{Q-m} (p)=\int_\LL \int_{R\cap u^{-1}(s)} h(p)\dd\psi^{Q-m-\ell}(p)\:\dd\psi^\ell(s) 
		\end{equation}
	holds.
\end{Corollary}

 \begin{Remark}
Let us stress that assumptions~\eqref{eq:technicalassumption} and~\eqref{eq:technicalassumptionPRODUCT} cannot be easily relaxed: given a map $u\in C^1_H(\Omega,\R^2)$ defined on an open subset $\Omega$ of the first Heisenberg group $\HH^1\equiv\R^3$, the validity of a coarea formula of the type
\[
\int_\Omega \calC(\DH u_p)\dd\psi^4(p) = \int_{\R^2}\psi^{2}(\Omega\cap u^{-1}(s))\dd\mathscr L^2(s)
\]
is indeed a challenging open problem as soon as $\DH u_p$ is surjective, see e.g.~\cite{Kozhevnikov,LeonardiMagnani,MagnaniStepanovTrevisan}. In our notation,  this situation corresponds to $\MM=\{0\}$ and $\LL=\R^2$. Since the kernel of any homogeneous surjective morphism $\HH^1\to\R^2$ is the center of $\HH^1$, which does not admit any complementary subgroup, 
 no point can be split-regular for $u$.
 Therefore, if~\eqref{eq:technicalassumptionPRODUCT} holds,
 then $\calC(\DH u_p)=0$ by Proposition~\ref{prop:coarea-factor},
 and thus both sides of the coarea formula are null.
 In particular,~\eqref{eq:technicalassumptionPRODUCT} implies that for $\cal L^2$-a.e.~$s\in\R^2$, $\psi^2(\Omega\cap u^{-1}(s))=0$. 
However, a coarea formula was proved for $u:\HH^n\to \R^{2n}$, assuming $u$ to be of class $C^{1,\alpha}_H$, see~\cite[Theorem~6.2.5]{Kozhevnikov} and also \cite[Theorem~8.2]{MagnaniStepanovTrevisan}.
\end{Remark}
\begin{Remark}
The following weak version of Sard's Theorem holds: under the assumptions and notation of Theorem~\ref{thm:coarea}, then
\begin{equation}\label{eq:sard1}
\psi^{Q-m-\ell}(\{p\in\Sigma:\DH (uf)_p(T^H_p\Sigma)\varsubsetneq\LL\}\cap u^{-1}(s))=0\quad \text{for $\psi^\ell$-a.e.~$s\in\LL$}.
\end{equation}
Moreover, since every level set $\Sigma\cap u^{-1}(s)$ is a $C^1_H$ submanifold around split-regular points of $uf$, Theorem~\ref{thm:coarea} implies that
\begin{equation}\label{eq:sard2}
\text{$\Sigma\cap u^{-1}(s)$ is $(\G;\KK)$-rectifiable\qquad for $\psi^\ell$-a.e.~$s\in\LL$.}
\end{equation}
Clearly, statements analogous to~\eqref{eq:sard1} and~\eqref{eq:sard2} hold under the assumptions and notation of either  Corollary~\ref{cor:coareaPRODUCT} or Theorem~\ref{thm:coareaHeisenberg} below.
\end{Remark}

The proof of Theorem~\ref{thm:coarea} follows the strategy used in~\cite{FedererCurvatureMeasures} (see also~\cite{MagnaniAreaCoarea}) and, as already mentioned, it stems from the area formula of Theorem~\ref{prop05161504}, as we now describe. First, in Proposition~\ref{prop07171752} we prove a coarea inequality, that in turn is based on an ``abstract'' coarea inequality (Lemma~\ref{lem-coarea-ineq}) for Lipschitz maps between metric spaces. Second, in Lemma~\ref{prop:coarea-factor} we prove Theorem~\ref{thm:coarea} in the ``linearized'' case when both $f$ and $u$ are homogeneous group morphisms: in this case formula~\eqref{eq:coarea2} holds with a constant coarea factor $\calC(\PP,L)$ which depends only on the normal homogeneous subgroup $\PP:=\ker f$ and on the homogeneous morphism $L=u$ (actually, on $L|_\Sigma$ only). Lemma~\ref{prop:coarea-factor}, whose proof  is a simple application of Theorem~\ref{prop05161504}, actually defines the coarea factor $\calC(\PP,L)$. The proof of Theorem~\ref{thm:coarea} is then a direct consequence of Theorem~\ref{thm:coarea2}, which states that for $\psi^{Q-m}$-a.e.~$p\in\Sigma$ the Federer density $\Theta_{\psi^d}(\mu_{\Sigma,u};p)$ of the measure
\[
\mu_{\Sigma,u}(E) \defeq \int_{\LL}  \psi^{Q-m-\ell}(E\cap\Sigma\cap u^{-1}(s)) \dd \cal\psi^\ell (s),\qquad E\subset\Omega
\]
is equal to $\calC(T^H_p\Sigma,\DH (uf)_p)$. For ``good'' points $p$, i.e., when $\DH  (uf)_p|_{T^H_p\Sigma}$ is  onto $\LL$, such equality is obtained by another application of  Theorem~\ref{prop05161504}, see Proposition~\ref{prop07251509}: this is the point where one needs the  assumption~\eqref{eq:technicalassumption}, which guarantees that, locally around good points, the level sets $\Sigma\cap u^{-1}(s)$ are $C^1_H$ submanifolds. The remaining  ``bad'' points, where $\DH  (uf)_p|_{T^H_p\Sigma}$ is not surjective on $\LL$, can be treated using the coarea inequality, see Lemma~\ref{lem07241228}.

Recall that the classical Euclidean coarea formula is proved  when the slicing function $u$ is only Lipschitz continuous.  Extending Theorem~\ref{thm:coarea} to the case where $u:\Sigma\to\LL$  is only Lipschitz seems for the moment out of reach. 
Observe that one should first provide, for a.e. $p\in\Sigma$, a notion of Pansu differential of $u$ on $T^H_p\Sigma$: this does not follow from Pansu’s Theorem~\cite{Pansu}.
Furthermore, the function $f$ in Theorem~\ref{thm:coarea} should play no role, and actually any result should depend only on the restriction of $u$ to $\Sigma$.

Let us also stress that, to the best of our knowledge, Theorem~\ref{thm:coarea} provides the first  sub-Riemannian coarea formula that is proved when the set $\Sigma$ is not a positive $\psi^Q$-measure subset of $\G$ (i.e., in the notation of Theorem~\ref{thm:coarea}, when $\MM=\{0\}$). The only exception to this is~\cite[Theorem~1.5]{MontiVittoneHeight}, where a coarea formula was proved for $C^1_H$ submanifolds of codimension 1 in Heisenberg groups $\HH^n,n\ge 2$. As a corollary of Theorem~\ref{thm:coarea}, we are able  both to extend this result, to all codimensions not greater than $n$, and to improve it, in the sense that we show that the implicit ``perimeter'' measures considered in~\cite[Theorem 1.5]{MontiVittoneHeight} on the level sets of $u$ are indeed Hausdorff or spherical measures. Furthermore, when $\HH^n$ is endowed with a rotationally invariant distance, $u$ takes values in $\R^\ell$, and the measures $\psi^d$ under consideration are $\cal S^d$, then the coarea factor coincides up to constants with the quantity
\begin{equation}\label{eq:defJRu}
 J^Ru(p):=\left( \det (L\circ L^T) \right)^{1/2},\qquad L:=\DH u_p|_{T^H_pR},
\end{equation}
In~\eqref{eq:defJRu}, the point $p$ belongs to a rectifiable set $R\subset\HH^n$ and, by abuse of notation, we use  standard exponential coordinates on $\HH^n\equiv\R^{2n+1}$ to identify $T^H_pR$ with a $(2n+1-m)$-dimensional plane; with this identification $\DH u_p$ is a linear map on $\R^{2n+1}$ that is, actually, independent of the last ``vertical'' coordinate. The superscript $^T$ denotes transposition. 

\begin{Theorem}[Coarea formula in Heisenberg groups]\label{thm:coareaHeisenberg}
	Consider  an open set $\Omega\subset\HH^n$, a $(\HH^n,\R^m)$-rectifiable set $R\subset\Omega$  and a function $u\in C^1_H(\Omega;\R^\ell)$ such that $1\le m+\ell\le n$.  
	Then, for every Borel function $h:R\to[0,+\infty)$ the equality
	\begin{equation*}
	\int_R h(p)\calC(T^H_pR,\DH u_p) \, \dd\psi^{2n+2-m} (p)=\int_{\R^\ell} \int_{R\cap u^{-1}(s)} h(p)\dd\psi^{2n+2-m-\ell}(p)\,\dd\psi^\ell(s) 
		\end{equation*}
	holds.
	
	Moreover, if $\HH^n$ is endowed with a rotationally invariant distance $\rho$, then there exists a constant $\textfrak c=\textfrak c(n,m,\ell,\rho)>0$ such that
		\begin{equation*}
	\textfrak c \int_R h(p)J^Ru(p) \, \dd\cal S^{2n+2-m} (p)=\int_{\R^\ell} \int_{R\cap u^{-1}(s)} h(p)\dd\cal S^{2n+2-m-\ell}(p)\,\dd\cal L^\ell(s).
		\end{equation*}
	\end{Theorem}
	The first  statement of Theorem~\ref{thm:coareaHeisenberg} is an immediate application of Corollary~\ref{cor:coareaPRODUCT}, while the second one needs an explicit representation for the spherical measure on {\it vertical} subgroups of $\HH^n$ (i.e., elements of $\scr T_{\HH^n,\R^k}$) which use results of~\cite{CorniMagnani}. See Proposition~\ref{prop:Sverticalsubgroups}.

\medskip

{\it Acknowledgements.} The authors are grateful to F.~Corni, V.~Magnani, R.~Monti and P.~Pansu for several stimulating discussions. They wish to express their gratitude to A.~Merlo for suggesting to address the density existence problem of Corollary~\ref{cor:densityexistence}.

\section{Preliminaries}\label{sec:preliminaries}

\subsection{First definitions}\label{sec12101536}
Let $V$ be a real vector space with finite dimension and  $[\cdot,\cdot]:V\times V\to V$ be the Lie bracket of a Lie algebra $\frk g=(V,[\cdot,\cdot])$.
We say that $\g$ is \emph{graded} if  subspaces $V_1,\dots,V_s$ are fixed so that
\begin{align*}
&V=V_1\oplus\dots\oplus V_s\\
\text{and }&[V_i,V_j]:=\Span\{[v,w]:v\in V_i,\ w\in V_j\}\subset V_{i+j}\text{ for all $i,j\in\{1,\dots,s\}$,}
\end{align*}
where we agree that $V_k=\{0\}$ if $k>s$. 
Graded Lie algebras are nilpotent.
A graded Lie algebra is \emph{stratified of step $s$} if equality  $[V_1,V_j]=V_{j+1}$ holds and $V_s\neq\{0\}$.
Our main object of study are stratified Lie algebras, but we will often work with subspaces that are only graded Lie algebras.

On the vector space $V$ we define a group operation via the Baker--Campbell--Hausdorff formula
\begin{align*}
pq &:=
\sum_{n=1}^\infty\frac{(-1)^{n-1}}{n} \sum_{\{s_j+r_j>0:j=1\dots n\}}
	\frac{ [p^{r_1}q^{s_1}p^{r_2}q^{s_2}  \cdots p^{r_n}q^{s_n}] }
	{\sum_{j=1}^n (r_j+s_j) \prod_{i=1}^n r_i!s_i!  } \\
	&=p+q+\frac12[p,q]+\dots ,
\end{align*}
where 
\[
[p^{r_1}q^{s_1}p^{r_2}q^{s_2}  \cdots p^{r_n}q^{s_n}]
= \underbrace{[p,[p,\dots,}_{r_1\text{ times}}
\underbrace{[q,[q,\dots,}_{s_1\text{ times}}
\underbrace{[p,\dots}_{\dots}]\dots]]\dots]] .
\]
The sum in the formula above is  finite  because $\frk g$ is nilpotent.
The resulting Lie group, which we denote by $\G$, is nilpotent and simply connected; we will call it \emph{graded group} or  \emph{stratified group}, depending on the type of grading of the Lie algebra.
The identification $\G=V=\frk g$ corresponds to the identification between Lie algebra and Lie group via the exponential map $\exp:\g\to\G$.
Notice that $p^{-1}=-p$ for every $p\in\G$ and that $0$ is the neutral element of $\G$.

If $\frk g'$ is another graded Lie algebra with underlying vector space $V'$ and Lie group $\G'$, then, with the same identifications as above, a map $V\to V'$ is a Lie algebra morphism if and only if it is a Lie group morphism, and all such maps are linear.
In particular, we denote by $\Hom_h(\G;\G')$ the space of all \emph{homogeneous morphisms} from $\G$ to $\G'$, that is, all linear maps $V\to V'$ that are Lie algebra morphisms (equivalently, Lie group morphisms) and that map $V_j$ to $V_j'$.
If $\g$ is stratified, then homogeneous morphisms are uniquely determined by their restriction to $V_1$.

For $\lambda>0$, define the \emph{dilations} as the maps $\delta_\lambda:V\to V$ such that $\delta_\lambda v=\lambda^j v$ for $v\in V_j$.
Notice that $\delta_\lambda\delta_\mu=\delta_{\lambda\mu}$
and that $\delta_\lambda\in\Hom_h(\G;\G)$, for all $\lambda,\mu>0$.
Notice also that a Lie group morphism $F:\G\to\G'$ is homogeneous if and only if $F\circ\delta_\lambda = \delta_\lambda'\circ F$ for all $\lambda>0$, where $\delta'_\lambda$ denotes the dilations in $\G'$.
We say that a subset $M$ of $V$ is \emph{homogeneous} if $\delta_\lambda(M)=M$ for all $\lambda>0$.
Let $\PP$ be a homogeneous subgroup of $\G$ and $\theta$ a Haar measure on $\PP$.
Since $\delta_\lambda|_{\PP}$ is an automorphism of $\PP$, there is $c_\lambda>0$ such that $(\delta_{\lambda})_\#\theta = c_\lambda\theta$.
Since the map $\lambda\mapsto \delta_\lambda|_{\PP}$ is a multiplicative one-parameter group of automorphisms, the map $\lambda\mapsto c_\lambda$ is a continuous automorphism of the multiplicative group $(0,+\infty)$, hence $c_\lambda = \lambda^{-d}$ for some $d\in\R$.
As $\delta_\lambda$ is contractive for $\lambda<1$, we actually have $d>0$.
Since any other Haar measure of $\PP$ is a positive multiple of $\theta$, the constant $d$ does not depend on the choice of the Haar measure.
We call such exponent $d$ the \emph{homogeneous dimension} of $\PP$.
The homogeneous dimension of the ambient space $\G$ is denoted by $Q$ and it is easy to see that $Q:=\sum_{i=1}^si\dim V_i$.

A \emph{homogeneous distance} on $\G$ is a distance function $\rho$ that is left-invariant and 1-homogeneous with respect to dilations, i.e., 
\begin{enumerate}
\item[(i)]
$\rho(gx,gy)=\rho(x,y)$ for all $g,x,y\in\G$;
\item[(ii)]
$\rho(\delta_\lambda x,\delta_\lambda y) = \lambda \rho(x,y)$ for all $x,y\in\G$ and all $\lambda>0$.
\end{enumerate}
When a stratified group $\G$ is endowed with a homogeneous distance $\rho$, we call the metric Lie group $(\G,\rho)$ a \emph{Carnot group}.
Homogeneous distances induce the  topology of  $\G$, see \cite[Proposition 2.26]{MR3943324}, and are biLipschitz equivalent to each other.
Every homogeneous distance defines a \emph{homogeneous norm}  $\|\cdot\|_\rho:\G\to[0,+\infty)$, $\|p\|_\rho := \rho(0,p)$.
We denote by $|\cdot|$ the Euclidean norm in $\R^\ell$. 
The following property relating norm and conjugation, proved in \cite[Lemma~2.13]{FranchiSerapioni2016Intrinsic}, will be useful: there exists $C=C(\G,\rho)>0$ such that
\begin{equation}\nonumber 
\|q^{-1}pq\|_\rho \le \|p\|_\rho+C\left( \|p\|_\rho^{1/s}\|q\|_\rho^{(s-1)/s}+\|p\|_\rho^{(s-1)/s}\|q\|_\rho^{1/s}\right)\quad\forall\: p,q\in\G.
\end{equation}

Open balls with respect to $\rho$ are denoted by $\Ball_\rho(x,r)$, closed balls by $ \cBall_\rho(x,r)$, or 
simply $\Ball(x,r)$ and $ \cBall(x,r)$
if it is clear which distance we are using.
We also use the notation $\cBall(E,r) := \{x: d(x,E)\le r\}$ for subsets $E$ of $\G$.
The diameter of a set with respect to $\rho$ is denoted by $\diam(E)$ 
or $\diam_\rho(E)$. 
Notice that $\diam_\rho(\Ball_\rho(p,r)) = 2r$, for all $p\in\G$ and $r>0$.  By left-invariance of $\rho$ it suffices to prove this for $p=0$.
On the one hand the triangle inequality implies $\diam_\rho(\Ball_\rho(0,r)) \le 2r$.
On the other hand, if $v\in V_1$ is such that $\rho(0,v)=r$, then $\rho(0,v^{-1})=r$ and $\rho(v^{-1},v) = \rho(0,2v) = 2\rho(0,v) = 2r$, because $vv = v+v = \delta_2v$.
It follows that $\diam_\rho(\Ball_\rho(0,r)) \ge 2r$. 

If $\rho$ and $\rho'$ are homogeneous distances on $\G$ and $\G'$, the distance between two homomorphisms $L,M\in\Hom_h(\G;\G')$ is
\[
d_{\rho,\rho'}(L,M) := \max_{p\neq0}\frac{\rho'(L(p),M(p))}{\|p\|_\rho} 
= \max_{\|p\|_\rho=1}\rho'(L(p),M(p)) .
\]
The function $d_{\rho,\rho'}$ is a distance on $\Hom_h(\G;\G')$ inducing the manifold topology.

\subsection{Measures and Federer density}\label{sec09172014}
In the following, the word measure will stand for outer measure. We work on $\G$ and its subsets endowed with the metric $\rho$. In particular, the balls are those defined by $\rho$ and the Hausdorff dimension of $(\G,\rho)$ coincides with the homogeneous dimension $Q$.

For $d\in [\,0,Q\,]$, let $\calH^d$ and $\calS^d$ be the Hausdorff and spherical Hausdorff measures of dimension $d$ in $\G$ defined for $E\subset\G$ by
\begin{align*}
&\calH^d(E):=\lim_{\epsilon\to 0^+}\quad\inf\left\{\sum_{j\in\N}(\diam E_j)^d:E\subset\bigcup_{j\in\N}E_j,\ \diam E_j<\epsilon\right\},\\
&\calS^d(E):=\lim_{\epsilon\to 0^+}\quad\inf\left\{\sum_{j\in\N}(2r_j)^d:E\subset\bigcup_{j\in\N}\cBall(x_j,r_j),\ 2r_j<\epsilon\right\}.
\end{align*}
It is clear that, in the definition of $\calH^d$, one can ask the covering sets $E_j$ to be closed.
Moreover, we clearly have  $\calH^d(E) \le \calS^d(E) \le 2^d \calH^d(E)$. 
Note that contrarily to the usual Euclidean or Riemannian definition, we do not introduce normalization constants; this is due to the fact that the appropriate constant is usually linked to the solution to the isodiametric problem, which is open in Carnot Groups and their subgroups and also highly dependent on the metric~$\rho$. See also Remark~\ref{rem:isodiametric}. In the following, $\psi^d$ will be either $\calH^d$ or $\calS^d$ and $\scr E$ will be, respectively,  the collection of closed subsets of $\G$ of positive diameter or the collection of closed balls in $\G$ with positive diameter.

If $\mu$ is a measure on $\G$, define the \emph{$\psi^d$-density of $\mu$ at $x\in\G$} as
\begin{align}
\Theta_{\psi^d}(\mu;x) 
&\defeq  \lim_{\epsilon\to 0^+}
\sup\left\{ \frac{\mu(E)}{(\diam E)^d} : x\in E\in\scr E,\ \diam E\le\epsilon \right\}.\label{eq:defFedererdensity}
\end{align}
This upper density is sometimes called \emph{Federer density} \cite{FSSCcentered,MagnaniEdinburgh,MagnaniNewDifferentiation}; note that if $\psi^d$ is the spherical measure, its Federer density can differ from the usual spherical density, as the latter involves centered balls. 
Recall that a measure $\nu$ is \emph{Borel regular} if open sets are measurable and 
for every $A\subset\G$ there exists a Borel set $A'\subset\G$ such that $A\subset A'$ and $\nu(A')=\nu(A)$.
We will use the following density estimates, which follow from {\cite[Theorems 2.10.17 and 2.10.18]{FedererGMT}}.

\begin{Theorem}[Density estimates]\label{thm:Federer}
  Let $\psi^d$ be as above,
 $\mu$ a Borel regular measure, and fix $t>0$ and a  set $A$ in $\G$. Then
  \begin{enumerate}
    \item[(i)] if $\Theta_{\psi^d}(\mu;x) <  t$  for all $x\in A$,  then $\mu(A) \le t \psi^d(A)$,
    \item[(ii)] if $\Theta_{\psi^d}(\mu;x)>  t$ for all $x\in A$,  then 
    $\mu(A)\ge t \psi^d(A)$.
  \end{enumerate}
\end{Theorem}
A direct consequence of these results is the following (see also \cite[Theorem~9]{MagnaniEdinburgh} and~\cite[Theorem 1.11]{FSSCcentered}).
\begin{Proposition}\label{propFedDensity}
	If $\mu$ is
	locally finite and Borel regular 
	on $\G$, and if $x\mapsto \Theta_{\psi^d}(\mu;x)$ is a Borel function which is positive and finite $\mu$-almost everywhere,
	then
	\[
	   \mu = \Theta_{\psi^d}(\mu;\cdot) \psi^d.
	\]
\end{Proposition}
Proving that the Federer density is a $\psi^d$-measurable or a Borel function is in general not an easy task; we  provide a criterion, which will be useful later in Sections~\ref{sec:goodpoints} and~\ref{sec:badpoints}.
Recall that a Borel measure $\nu$  is \emph{doubling} if there exists $C\ge 1$ such that $\nu(\Ball(p,2r))\le C\,\nu(\Ball(p,r))$ for all $p\in \G$ and $r>0$.

\begin{Proposition}\label{prop:RadonNykodym}
Given a set $\Sigma\subset \G$ such that $\psi^d\hel \Sigma$ is locally doubling Borel regular measure,
 assume that $\mu$ is a locally finite Borel regular measure, absolutely continuous with respect to $\psi^d\hel \Sigma$; then
\begin{enumerate}[label=(\roman*)]
\item $\Theta_{\psi^d}(\mu;\cdot)$ is $(\psi^d\hel\Sigma)$-measurable;
\item $\Theta_{\psi^d}(\mu;\cdot)<+\infty$, $\psi^d$-a.e.~on $\Sigma$ and
\[
\Theta_{\psi^d}(\mu;p)=\lim_{r\to 0^+}\frac{\mu(\cBall(p,r))}{\psi^d(\Sigma\cap\cBall(p,r))},\quad \text{for $\psi^d$-a.e.~}p\in\Sigma;
\]
\item $\mu=\Theta_{\psi^d}(\mu;\cdot) \psi^d\hel\Sigma$.
\end{enumerate}
In particular
\begin{equation}\nonumber
\lim_{r\to 0^+}\aveint{\Sigma\cap\cBall(p,r)}{}\left|\Theta_{\psi^d}(\mu;\cdot) - \Theta_{\psi^d}(\mu;p)\right|\dd \psi^d=0,\qquad\text{for $\psi^d$-a.e.~}p\in\Sigma. 
\end{equation}
\end{Proposition}
\begin{proof}
It is well-known (see e.g.~\cite{rigot2018differentiation}) that Radon-Nikodym Differentiation Theorem holds for differentiating a measure with respect to a doubling measure. 
Precisely, by combining~\cite[Theorems 2.2,~2.3,~3.1]{rigot2018differentiation} one infers that the Radon-Nikodym derivative
\[
\Theta(p):=\lim_{r\to 0^+}\frac{\mu(\cBall(p,r))}{\psi^d(\Sigma\cap\cBall(p,r))}
\]
exists and is finite $\psi^d$-a.e.~on $\Sigma$. Moreover, $\Theta$ is $(\psi^d\hel\Sigma)$-measurable,  $\mu=\Theta \psi^d\hel\Sigma$ and (see~\cite[Section 2.7]{Heinonen})
\[
\lim_{r\to o^+}\aveint{\Sigma\cap\cBall(p,r)}{}\left|\Theta - \Theta(p)\right|\dd \psi^d=0\qquad\text{for $\psi^d$-a.e.~}p\in\Sigma. 
\]
As a consequence, we have only to prove that $\Theta_{\psi^d}(\mu;p)=\Theta(p)$ for $\psi^d$-a.e.~$p\in \Sigma$. In turn, it is enough to show that, for every fixed $s,t\in\Q$, $s<t$, the sets
\begin{align*}
& A:=\{p\in\Sigma:\Theta(p)<s<t<\Theta_{\psi^d}(\mu;p)\}\\
& B:=\{p\in\Sigma:\Theta_{\psi^d}(\mu;p)<s<t<\Theta(p)\}
\end{align*}
are $\psi^d$-negligible. 
On the one hand, let $A'$ be a Borel set with $A\subset A'$, $\psi^d(A)=\psi^d(A')$ and $A'\subset\{\Theta<s\}$.
Then
\begin{align*}
s\psi^d(A)
= s\psi^d(A')
\ge \int_{A'} \Theta\dd\psi^d
= \mu(A')
\ge \mu(A)\ge  t\psi^d(A),
\end{align*}
where the last inequality is a consequence of Theorem~\ref{thm:Federer} $(ii)$.
Therefore, $\psi^d(A)=0$. On the other hand, let $B'$ be a Borel set with $B\subset B'$, $\mu(B')=\mu(B)$. 
Then
\[
t\psi^d(B) 
\le \int_{B'}\Theta\dd\psi^d
= \mu(B')
= \mu(B) 
\le s\psi^d(B),
\]
where the last inequality is a consequence of Theorem~\ref{thm:Federer} $(i)$.
Therefore $\psi^d(B)=0$.
\end{proof}

\subsection{Pansu differential}\label{sec:Pansudifferentiability}
Let $\G$ and $\G'$ be two Carnot groups and $\Omega\subset\G$ open.
A function $f:\Omega\to\G'$ is \emph{Pansu differentiable at $p\in\Omega$} if 
there is $L\in\Hom_h(\G;\G')$ 
such that
\[
\lim_{x\to p}
	\frac{\rho'( f(p)^{-1}f(x) , L(p^{-1}x) )}{\rho(p,x)}
= 0 .
\]
The map $L$ is called \emph{Pansu differential} of $f$ at $p$ and it is denoted by $\DH f(p)$ or $\DH f_p$.
A map $f:\Omega\to\G'$ is  \emph{of class $C^1_H$}
if $f$ is Pansu differentiable at all points of $\Omega$ and the Pansu differential $p\mapsto \DH f(p)$ is continuous.
We denote by $C^1_H(\Omega;\G')$ the space of all  maps from $\Omega$ to $\G'$ of class $C^1_H$.

A function $f:\Omega\to\G'$ is \emph{strictly Pansu differentiable at $p\in\Omega$} if 
there is $L\in\Hom_h(\G;\G')$ 
such that
\[
\lim_{\epsilon\to0}\ \sup \left\{
	\frac{\rho'( f(y)^{-1}f(x) , L(y^{-1}x) )}{\rho(x,y)}
	: x,y\in\Ball_\rho(p,\epsilon),\ x\neq y \right\}
= 0 .
\]
Clearly, in this case $f$ is Pansu differentiable at $p$ and $L=\DH f(p)$.

The next results allows us to simplify several arguments in the sequel:
\begin{Proposition}\label{prop:C1Hiffstrict}
A function $f:\Omega\to\G'$ is of class $C^1_H$ on $\Omega$ if and only if $f$ is strictly Pansu differentiable at all points in $\Omega$.
\end{Proposition}
\begin{proof}
Assume that $f\in C^1_H(\Omega,\G')$ and let $p\in\Omega$ be fixed; then, by \cite[Theorem 1.2]{MR3123745} one has
\[
\lim_{\epsilon\to0}\ \sup \left\{
	\frac{\rho'( f(y)^{-1}f(x) , \DH f_x(y^{-1}x) )}{\rho(x,y)}
	: x,y\in\Ball_\rho(p,\epsilon),\ x\neq y \right\}=0\,.
\]
The continuity of $x\mapsto \DH f_x$ provides
\[
\lim_{\epsilon\to0}\ \sup \left\{
	\frac{\rho'(  \DH f_x(y^{-1}x) , \DH f_p(y^{-1}x) )}{\rho(x,y)}
	: x,y\in\Ball_\rho(p,\epsilon),\ x\neq y \right\}=0
\]
and the strict differentiability of $f$ at $p$ follows.

Conversely, assume that $f$ is strictly Pansu differentiable at all points in $\Omega$; we have to prove that $p\mapsto \DH f_p$ is continuous. Assume not, i.e., assume there exist $\delta>0$ and, for every $n\in\N$, points $x_n\in\Omega$ and $v_n\in\G$ such that $\|v_n\|_\rho=1$, $x_n\to p$ and
\[
\rho'(\DH f_{x_n}(v_n),\DH  f_p(v_n))\ge 2\delta\qquad\forall\:n\in\N.
\]
By strict differentiability of $f$ at $p$ there exist $\bar n$ and $\bar s>0$ such that
\[
\frac{\rho'(f(x_n)^{-1}f(x_n\delta_s v_n),\DH  f_p(\delta_s v_n))}{s}\le \delta\qquad\forall\: n\ge\bar n,\ s\in(0,\bar s).
\]
In particular, for every $n\ge\bar n$ and $s\in(0,\bar s)$ we have
\begin{eqnarray*}
 {\rho'(f(x_n)^{-1}f(x_n\delta_s v_n),\DH  f_{x_n}(\delta_s v_n))}\hspace{-4.5cm}&\\
&\ge&\,  {\rho'(\DH  f_{p}(\delta_s v_n),\DH  f_{x_n}(\delta_s v_n))}
-{\rho'(f(x_n)^{-1}f(x_n\delta_s v_n),\DH  f_{p}(\delta_s v_n))}\\
&\ge&\,  2\delta s-\delta s=\delta s.
\end{eqnarray*}
This would contradict the differentiability of $f$ at $x_n$.
\end{proof}

\begin{Lemma}\label{lem10051109}
	If $f\in C^1_H(\Omega;\G')$, then $f:(\Omega,\rho)\to(\G',\rho')$ is locally Lipschitz.
\end{Lemma}
\begin{proof}
	Let $p\in\Omega$.
	By strict differentiability of $f$ at $p$, there is $\epsilon>0$ such that 
	\[
	\frac{\rho'( f(y)^{-1}f(x) , L(y^{-1}x) )}{\rho(y,x)}<1 \qquad\text{for all $x,y\in \Ball_\rho(p,\epsilon)$, $x\neq y$,}
	\]
	 where $L=\DH f(p)$.
	Since $\rho'(0,L  (y^{-1}x))\le C \rho(y,x)$ for some positive $C$, 
	then $\rho'(f(y),f(x))=\rho'(0,f(y)^{-1}f(x)) \le (C+1)\rho(y,x)$, that is, $f$ is Lipschitz continuous on $\Ball_\rho(p,\epsilon)$.
	\end{proof}

\begin{Remark}\label{rem:rem2.3}
The notion of Pansu differentiability, as well as Lemma \ref{lem10051109}, can be stated also when the target group $\G'$ is only graded. However, there is no loss of generality in assuming $\G'$ to be stratified. Indeed, if $f:\Omega\to\G'$ is locally Lipschitz, then the image of a rectifiable curve in $\G$ is a rectifiable curve in $\G'$ tangent to the first layer $V_1'$ in the grading of $\G'$; since $\G$ is stratified, each connected component $U$ of $\Omega$ is pathwise connected by rectifiable curves, and this implies that $f(U)$ is contained in (a coset of) the stratified subgroup of $\G'$ generated by $V_1'$.

Moreover, as soon as $f$ is  open, or has a regular point, then $\G'$ must be a Carnot group.
\end{Remark}

\subsection{Intrinsic graphs and implicit function theorem}\label{sec09181614}
We refer to \cite{FranchiSerapioni2016Intrinsic} for a more general theory of intrinsic graphs.
Recall the identification $\G=\frk g=V$ that we  made in Section~\ref{sec12101536}.

\begin{Lemma}\label{lem09151008}
	Let $\V$ and $\W$ be homogeneous linear subspaces of a graded group $\G$.
	If $\V\cap\W=\{0\}$ and $\dim\V+\dim\W=\dim\G$, then the map $\W\times\V\to\G$, $(w,v)\mapsto wv$, is a surjective diffeomorphism.
\end{Lemma}
\begin{proof}
	Denote by $\phi:\W\times\V\to\G$ the map $\phi(w,v):=wv$.
	Since its differential at $(0,0)$ is a linear isomorphism, $\phi$ is a diffeomorphism from a neighborhood of $(0,0)$ to a neighborhood of $0\in\G$.
	Since $\phi(\delta_\lambda w,\delta_\lambda v) = \delta_\lambda\phi(w,v)$ for all $\lambda>0$, we conclude that $\phi$ is a surjective diffeomorphism onto $\G$.
\end{proof}

A homogeneous subgroup $\W$ is \emph{complementary} to a homogeneous subgroup $\V$ if $\G=\W\V$ and $\W\cap\V=\{0\}$.
We denote by $\mathscr W_\V$ the set of all homogeneous subgroups of $\G$ that are complementary to $\V$.
By Lemma~\ref{lem09151008}, we have $\W\in\mathscr W_\V$ if and only if $\V\in\mathscr W_\W$.
Again by Lemma~\ref{lem09151008}, any choice of $\V$ and $\W\in\mathscr W_\V$ gives two projections
\begin{equation}\label{eq:defproiezioni}
\pi_\W:\G\to\W,\qquad \pi_\V:\G\to\V,
\end{equation}
which are defined, for every $p\in\G$, by requiring  $\pi_\W(p)=w\in\W$ and $\pi_\V(p)=v\in\V$ to be the only elements such that $p=wv$.
We will also write $p_\W$ and $p_\V$ for $\pi_\W(p)$ and $\pi_\V(p)$, respectively.

We say that a normal homogeneous subgroup \emph{$\W$ splits $\G$} if $\scr W_\W\neq\emptyset$.
In this case we call a choice of $\W$ and $\V\in\scr W_\W$ a \emph{splitting} of $\G$ and we write $\G=\W\cdot\V$. 
We say that $p\in\Omega$ is a \emph{split-regular point} of $f$ if the Pansu differential of $f$ at $p$ exists and is surjective, and if $\ker(\DH f(p))$ splits $\G$.
Recall that the kernel of a group morphism is always normal.
A \emph{singular point} is a point that is not split-regular.

\begin{Remark}\label{rem:V=G'} 
We observe that, if $p\in\Omega$ is a split-regular point of $f\in C^1_H(\Omega;\G')$ and $\V\in\mathscr W_{\ker(\DH f(p))}$,
then $\DH f(p)|_{\V}:\V\to\G'$ is an isomorphism of graded groups.
In particular, $\V$ is necessarily stratified.
For instance, if $\G'=\R^m$, then $\V$ is an Abelian subgroup of $\G$ contained in $V_1$.
\end{Remark}

 Notice  that a point can fail to be split-regular for $f\in C^1_H(\Omega;\G')$ for two distinct reasons: non-surjectivity of the differential, or non-existence of a splitting of $\G$ with the kernel of $\DH  f_p$ at some point $p$. However, the set of split-regular points  is open, i.e., if $\DH f_p$ is surjective and $(\ker\DH f_p)\cdot\V$ is a splitting, then, for $q$ close enough to $p$,  $\DH f_q$ is surjective and $(\ker\DH f_q)\cdot\V$ is  a splitting.

\begin{Lemma}[Coercivity]\label{lem08011611}
	If $f\in C^1_H(\Omega;\G')$, $p\in\Omega$ is a split-regular point and $\V$ is complementary to $\ker(\DH f(p))$, then there are a neighborhood $U$ of $p$ and $C>0$ such that, for all $q\in U$ and $v\in\V$ with $qv\in U$, 
	\[
	\rho'(f(q),f(qv)) \ge C \|v\|_\rho .
	\]
\end{Lemma}
\begin{proof}
	Arguing by contradiction, assume that there are sequences $q_j\in\Omega$ and $v_j\in\V\setminus\{0\}$ such that $q_j\to p$ and $v_j\to 0$ as $j\to\infty$,
	and $\rho'(f(q_j),f(q_jv_j)) \le  \|v_j\|_\rho/j$.
	Up to passing to a subsequence, we can assume that there exists
	$\bar w = \lim_{j\to\infty} \delta_{\|v_j\|^{-1}} v_j$.
	It follows that $\bar w\in\V$ and $\|\bar w\|_\rho =1$.
	Moreover, by strict differentiability
	\[
	\DH f(p)\bar w = \lim_{j\to\infty} \frac{f(q_j)^{-1}f(q_jv_j)}{\|v_j\|_\rho} = 0 ,
	\]
	in contradiction with the fact that $\V$ is complementary to the kernel of $\DH f(p)$.
\end{proof}

Let $\W\in\mathscr W_{\V}$.
A set $\Sigma\subset\G$ is \emph{an intrinsic graph $\W\to\V$} if there is a subset $A\subset\W$ and a function $\phi:A\to\V$ such that $\Sigma=\{w\phi(w):w\in A\}$.
Clearly, $\Sigma\subset\G$ is an intrinsic graph $\W\to\V$ if and only if the map $\pi_\W|_\Sigma:\Sigma\to\W$ is injective; in particular, every $\PP\in\mathscr W_\V$ is an intrinsic graph $\W\to\V$.
Left translations and dilations of $\W\to\V$ intrinsic graphs are again $\W\to\V$ intrinsic graphs, see \cite[Proposition 3.6]{ArenaSerapioni}.

The proof of the following lemma is inspired by \cite[Theorem A.5]{MR3906289}. Similar statements are contained in  \cite[Theorem~3.27]{FSSCAdvMath} and \cite[Theorem 1.4]{MR3123745}.

\begin{Lemma}[Implicit Function Theorem]\label{lem05161502}\label{lem:implicitfunction}
	Let $\Omega_0\subset\G$ be open, $g\in C^1_H(\Omega_0;\G')$ and let $o\in\G$ be a split-regular point of $g$.
	Let $\G=\W\cdot\V$ be a splitting of $\G$ such that $\ker(\DH g(o))$ is an intrinsic graph $\W\to\V$. 
	Then there are neighborhoods $A$ of $\pi_\W(o)$ in $\W$, $B$ of $g(o)$ in $\G'$ and $\Omega\subset\Omega_0$ of $o$, and a map $\varphi:A\times B\to\V$ such that the map $(a,b)\mapsto a\varphi(a,b)$ is a homeomorphism $A\times B\to \Omega$ and $g(a\varphi(a,b)) = b$.
	In particular,the map $\phi:A\to\V$ defined by $\phi(a):=\varphi(a,g(o))$ is such that $\{p\in\Omega:g(p)=g(o)\}=\{a\phi(a)\in\G:a\in A\}$.
\end{Lemma}

\begin{Remark}
 Notice that $\V\cap \ker(\DH g(o))=\{0\}$; hence, in the above lemma one can of course choose $\W=\ker(\DH g(o))$.
\end{Remark}

\begin{proof}[Proof of Lemma~\ref{lem05161502}]
	First, we prove that there is an open neighborhood $U\subset\Omega_0$ of $o$ such that the restriction $g|_{p\V}:p\V\cap U\to\G'$ is injective, for all $p\in U$.
	Arguing by contradiction, suppose that this is not the case.
	Then there are sequences $p_j,q_j\to o$ such that $p_j^{-1}q_j\in\V$ and $g(p_j)=g(q_j)$.
	From the strict Pansu differentiability of $g$ at $o$, it follows that
	\[
	0 = \lim_{j\to\infty} \frac{\rho'( g(q_j)^{-1}g(p_j) , \DH g(o)[p_j^{-1}q_j] ) }{ \rho(p_j,q_j) }
	= \lim_{j\to\infty}\left\| \DH g(o) \left[\delta_{\frac1{\rho(p_j,q_j)}} (p_j^{-1}q_j) \right] \right\|_{\rho'} .
	\]
	By the compactness of the sphere $\{v\in\V:\|v\|_\rho=1\}$, up to passing to a subsequence, there is $v\in\V$ with $\|v\|_\rho=1$ such that $\lim_{j\to\infty}\delta_{\rho(p_j,q_j)^{-1}} (p_j^{-1}q_j) = v$.
	Therefore, we obtain $\DH g(o)v=0$, in contradiction with the assumptions.
	This proves the first claim.
	
	Second, since the restriction $g|_{p\V\cap U}:p\V\cap U\to \G'$ is a continuous and injective map, and since both $\V$ and $\G'$ are topological manifolds of the same dimension, then we can apply the Invariance of Domain Theorem and obtain that $g|_{p\V\cap U}:p\V\cap U\to g(p\V\cap U)$ is a homeomorphism and that $g(p\V\cap U)$ is an open set.
		
	Third, let $U_2\Subset U_1\Subset U$ be open neighborhoods of $o$.
	We claim that there is $A\subset\W$ open such that $\pi_\W(o)\in A$ and such that for every $p\in o\V\cap U_2$ and for every $a\in A$ there is $q\in a\V\cap U_1$ such that $g(p)=g(q)$.
	Arguing by contradiction, suppose that this is not the case.
	Then there are sequences $a_j\in\W$ with $a_j\to\pi_\W(o)$ and $p_j\in o\V\cap U_2$ such that $g(p_j)\notin g(a_j\V\cap U_1)$.
	By the compactness of $\bar U_2$ and the continuity of $g$, for each $j$ there is $q_j\in a_j\V\cap \bar U_1$ such that
	\begin{equation}\label{eq12101651}
	\rho'(g(p_j),g(q_j)) = \inf\{\rho'(g(p_j),g(q)) : q\in a_j\V\cap U_1 \} .
	\end{equation}
	Since $g$ is a homeomorphism on each fiber $p\V\cap U$ and since $g(p_j)\notin g(a_j\V\cap U_1)$, we have $q_j\in a_j\V\cap\partial U_1$.
	Up to passing to a subsequence, there are $p_0\in o\V\cap \bar U_2$ and $q_0\in o\V\cap\partial U_1$ such that $p_j\to p_0$ and $q_j\to q_0$.
	Now, notice that $a_j\pi_\W(o)^{-1}\to 0$ and that, for $j$ large enough, we have $a_j\pi_\W(o)^{-1}p_j\in a_j\V\cap U_1$.
	Therefore, using~\eqref{eq12101651},
	\[
	\lim_{j\to\infty} \rho'( g(p_j) , g(q_j) )
	\le \lim_{j\to\infty} \rho'( g(p_j) ,g(a_j\pi_\W(o)^{-1}p_j) )
	=0 ,
	\]
	that is, $g(p_0)=g(q_0)$.
	Since $p_0\in  o\V\cap\bar U_2$ and $q_0\in o\V\cap (U\setminus U_1)$, this contradicts the injectivity of $g$ on $o\V\cap U$ and proves the  claim.
	
	Next, let $B:=g(o\V\cap U_2)$, which is an open neighborhood of $g(o)$, and $\Omega := \pi_\W^{-1}(A)\cap g^{-1}(B) \cap U_1$.
	The  previous claims imply that for every $a\in A$ and every $b\in B$ there is a unique $v\in\V$ such that $av\in \Omega$ and $g(av)=b$.
	Define $\varphi:A\times B\to \V$ as $\varphi(a,b)=v$.
	
	Finally, we claim that the map $\Phi(a,b):= a\varphi(a,b)$ is a homeomorphism $A\times B\to \Omega$.
	Notice that, if $p=\Phi(a,b)$, then $a=\pi_\W(p)$ and $b=g(p)$: therefore, $\Phi$ is injective.
	Moreover, if $p\in \Omega$, then $\pi_\W(p)\in A$, $g(p)\in B$ and $\Phi(\pi_\W(p),g(p)) = p$: therefore, $\Phi$ is also surjective.
	Finally, since $\Phi^{-1}:\Omega\to A\times B$ is a continuous bijection, then it is a homeomorphism by the Invariance of Domain Theorem.
	This completes the proof.
\end{proof}

We observe that, when $g:\G\to\G'$ is a homogeneous group morphism, then the statement of Lemma~\ref{lem05161502} holds with $A=\W$, $B=\G'$ and $\Omega=\G$.

\begin{Lemma}\label{lem09122015}
	Under the assumptions and notation of Lemma~\ref{lem05161502},
	suppose $o=0$ and
	define for $\lambda>0$
	\[
	\begin{array}{lccl}
	\varphi_{\lambda}:&\delta_{1/\lambda}A \times\delta_{1/\lambda}B &\to &\delta_{1/\lambda}\Omega \\
	& (a,b) &\mapsto & \delta_{1/\lambda}\varphi(\delta_\lambda a,\delta_\lambda b)
	\end{array}
	\]
	Let $\varphi_0$ be the implicit function associated with $\DH  g(0)$, that is, $\varphi_0:\mathbb W\times\G'\to\mathbb V$ is such that $\DH g(o)(a\varphi_0(a,b))=b$ for all $a$ and $b$.
	
	Then $\varphi_{\lambda}\to \varphi_0$ locally uniformly as $\lambda\to0^+$.
\end{Lemma}

\begin{proof}
	Without loss of generality, we assume $\Omega$ to be compactly contained in the domain of $g$.
	Define $g_\lambda:\delta_{1/\lambda}\Omega\to\G'$ by 
	\[
	g_\lambda(x) = \delta_{1/\lambda}g(\delta_\lambda x) .
	\]
	Notice that $g_\lambda(a\varphi_\lambda(a,b)) = b$ for all $(a,b)\in\delta_{1/\lambda}A \times\delta_{1/\lambda}B$.
	Possibly taking a smaller $\Omega$, by Lemma~\ref{lem08011611} there is $C>0$ such that $\rho'(g(x),g(y))\ge C\rho(x,y)$ for all $x,y\in\Omega$ with $\pi_\W(x)=\pi_\W(y)$.
	It follows that that $\rho'(g_\lambda(x),g_\lambda(y)) \ge C\rho(x,y)$ for all $x,y\in\delta_{1/\lambda}\Omega$ with $\pi_\W(x)=\pi_\W(y)$, because $\pi_\W\circ\delta_\lambda=\delta_\lambda\circ\pi_\W$.
	
	Fix a compact set $K\subset\W\times\G'$ and let $(a,b)\in K$. Then, for large enough $\lambda$ (depending only on $K$) we have $(a,b)\in\delta_{1/\lambda}A \times\delta_{1/\lambda}B$, $a\varphi_\lambda(a,b)\in\delta_{1/\lambda}\Omega$ and $a\varphi_0(a,b)\in\delta_{1/\lambda}\Omega$, hence
	\begin{align*}
	\rho(\varphi_\lambda(a,b),\varphi_0(a,b))
	&= \rho(a\varphi_\lambda(a,b),a\varphi_0(a,b)) \\
	&\le\frac1C \rho'( g_\lambda(a\varphi_\lambda(a,b)) , g_\lambda(a\varphi_0(a,b)) ) \\
	&=\frac1C \rho'( b , g_\lambda(a\varphi_0(a,b)) ) \\
	&=\frac1C \rho'( \DH g_0(a\varphi_0(a,b)) , g_\lambda(a\varphi_0(a,b)) ) .
	\end{align*}
	Since $g$ is Pansu differentiable at 0, $g_\lambda\to \DH g_0$ uniformly on compact sets.
	The map $(a,b)\mapsto a\varphi_0(a,b)$ is a homeomorphism $\W\times\G'\to\G$,
	hence $\varphi_\lambda\to\varphi_0$ uniformly on compact sets.
\end{proof}

\subsection{\texorpdfstring{$C^1_H$}{C1H} submanifolds and rectifiable sets}
\label{sec:C1H_and_rectifiable}
A set $\Sigma\subset\G$ is a \emph{submanifold of class }$C^1_H$ (or $C^1_H$ submanifold for short)  if there exists a Carnot group $\G'$ such that for every $p\in\Sigma$ there are an open neighborhood $\Omega$ of $p$ in $\G$ and a function $f\in C^1_H(\Omega;\G')$ such that 
$p$ is split-regular for $f$ and $\Sigma\cap\Omega = \{f=0\}$.
In this case, we sometimes call $\Sigma$ a $C^1_H(\G;\G')$-submanifold.

The \emph{homogeneous tangent subgroup} to $\Sigma$ at $p\in\Sigma$ is the homogeneous normal subgroup $T^H_p\Sigma := \ker(\DH f(p))$.
Statement $(iii)$ in the next lemma implies that $T^H_p\Sigma$ does not depend on the choice of $f$. Observe also that the homogeneous dimension of $T_p^H\Sigma$ is equal to the difference of the homogeneous dimensions of $\G$ and $\G'$ and is, in particular, independent of $p$; we call this integer \emph{homogeneous dimension} of $\Sigma$ and denote it by $\dim_H\Sigma$.

\begin{Definition}\label{def:C1WV}
Given a splitting $\G=\W\cdot\V$ and an open set $A\subset\W$, we say that \emph{$\phi:A\to\V$ is of class $C^1_{\W,\V}(A)$} if the intrinsic graph $\Sigma$ of $\phi$ is a $C^1_H$ submanifold and $T^H_{w\phi(w)}\Sigma\in \mathscr W_\V$ for every $w\in A$. 
\end{Definition}

Observe that, since $\V$ is isomorphic to $\G'$, the homogeneous dimension of $\W$ is equal to that of $\Sigma$.

\begin{Lemma}\label{lem09151054}
	Let $\Sigma\subset\G$ be a $C^1_H$ submanifold and $o\in\Sigma$.
	Let $\G=\W\cdot\V$ be a splitting such that $T^H_o\Sigma$ is the intrinsic graph of $\phi_0:\W\to\V$.
	The following statements hold:
	\begin{enumerate}[label=(\roman*)]
	\item
	There are open neighborhoods $\Omega$ of $o$ and $A$ of $\pi_\W(o)$, and a function $\phi\in C^1_{\W,\V}(A)$ such that $\Sigma\cap\Omega$ is the intrinsic graph of $\phi$.
	\item
	Define $\phi_\lambda(x) := \delta_{1/\lambda}\phi(\delta_\lambda x)$;
	then $\phi_\lambda\in C^1_{\W,\V}(\delta_{1/\lambda}A)$ and $\phi_\lambda\to\phi_0$ uniformly on compact sets as $\lambda\to0^+$.
	\item
	$\lim_{\lambda\to0^+} \delta_{1/\lambda}(o^{-1}\Sigma) = T^H_o\Sigma$ in the sense of local Hausdorff convergence of sets.
	\item
	If $U$ is a neighborhood of $o$ such that $\Sigma\cap U$ is the level set of $f\in C^1_H(U,\G')$ and $o$ is a split-regular point of $f$, then $\G'$ is isomorphic to $\V$.
	\end{enumerate}
\end{Lemma}
The proof of statements $(i)$, $(ii)$ and $(iii)$ is left to the reader, since it is a  consequence of Lemmas~\ref{lem:implicitfunction} and~\ref{lem09122015}. As for statement $(iv)$, it is enough to notice that the group morphism $\DH f(o)|_{\V}:\V\to\G'$ is injective (because $\V\cap\ker \DH f(o)=\{0\}$) and surjective (because $o$ is split-regular). 

An important property of the parametrizing map $\phi$ is that it is intrinsic Lipschitz in accordance with the theory developed by B.~Franchi, R.~Serapioni and F.~Serra Cassano, see e.g.~\cite{FrSeSe2006intrinsic,FranchiSerapioni2016Intrinsic}. 
We recall  that, given a splitting $\G=\W\cdot\V$ and $A\subset\W$, a map $\phi:A\to\V$ is \emph{intrinsic Lipschitz} if there exists $\cal C\subset \G$ such that the following conditions hold
\begin{itemize}
\item[(a)] $\cal C$ is a cone i.e., $\delta_\lambda\cal  C=\cal C$ for all $\lambda\ge 0$;
\item[(b)] $\V$ is an axis of $\cal C$, i.e., $\V\subset\cal  C$ and  $\V\setminus\{0\}\subset\mathring{\cal  C}$;
\item[(c)] the graph $\Sigma:=\{a\phi(a):a\in A\}$ of $\phi$ satisfies $\Sigma\cap (p\cal  C)=\{p\}$ for all $p\in\Sigma$.
\end{itemize}

\begin{Remark}
The above definition of  Lipschitz continuity for intrinsic graphs $\W\to\V$, tough slightly different, is equivalent to the one introduced by B.~Franchi, R.~Serapioni and F.~Serra Cassano, see e.g.~\cite[Remark A.2]{MR3906289}.
\end{Remark}

\begin{Corollary}\label{cor:intrLip}
	Intrinsic $C^1_H$ submanifolds are locally intrinsic Lipschitz graphs.
\end{Corollary}
\begin{proof}
	Let $\Sigma\subset\G$ be a $C^1_H$ submanifold, $o\in\Sigma$ and $\V\in\scr W_{T^H_o\Sigma}$.
	We need to prove that
	then there are a neighborhood $\Omega$ of $o$ and a cone $\cal C$ with axis $\V$ such that for all $p\in\Sigma\cap\Omega$ we have $(\Sigma\cap\Omega)\cap\cal C = \{p\}$.
	
	Let $\Omega$ be a neighborhood of $o$ with $f\in C^1_H(\Omega;\G')$ such that $\Sigma\cap \Omega=\{p\in\Omega:f(p)=f(o)\}$ and all points in $\Omega$ are split-regular for $f$.
	Up to shrinking $\Omega$, we can also assume, by Lemma~\ref{lem08011611}, that there exists $C>0$ such that
	\[
	\rho'(f(p),f(pv))\ge C\|v\|_\rho\quad\forall\: p\in\Omega,\ v\in\V\text{ such that }pv\in\Omega,
	\]
	and that, by Lemma~\ref{lem10051109}, $f:(\Omega,\rho)\to(\G',\rho')$ is $L$-Lipschitz, for some $L\ge0$.
	Define the cone
	\[
	\cal C:=\{0\}\cup\bigcup_{v\in\V} \Ball_\rho(v,\tfrac CL\|v\|_\rho)\subset \G.
	\]
	Requirements (a) and (b) above are clearly satisfied;
	to prove (c), let $p\in\Sigma\cap\Omega$ and $q\in(\Sigma\cap\Omega)\cap\cal C$.
	Then there exists $v\in \V$ such that $\rho(q,pv)<\tfrac CL\|v\|_\rho$, hence
	\[
	\rho'(f(p),f(q))
	\ge \rho'(f(p),f(pv))-\rho'(f(pv),f(q))\ge C\|v\|_\rho-L\rho(q,pv)>0.
	\]
	We conclude that $f(q)\neq f(p)$ and thus $q\notin\Sigma$. This completes the proof.
\end{proof}

The following result is an easy consequence of Lemma~\ref{lem09151054}, Corollary~\ref{cor:intrLip} and~\cite[Theorem~3.9]{FranchiSerapioni2016Intrinsic}. We denote by $\psi^d$ either the $d$-dimensional Hausdorff or $d$-dimensional spherical Hausdorff measure on $\G$ as in Section~\ref{sec09172014}.

\begin{Proposition}[Local Ahlfors regularity of the surface measure on $C^1_H$ submanifolds]\label{prop:doubling}
Let $\Sigma\subset\G$ be a $C^1_H$ and let $d:=\dim_H\Sigma$; then, for every compact set $K\subset\Sigma$ there exists $C=C(K)>0$ such that
\begin{equation}\label{eq:Ahlfors}
\frac1C r^d\le\psi^d(\Sigma\cap\Ball(p,r))\le Cr^d\qquad\forall\: p\in K.
\end{equation}
In particular, the measure $\psi^d\hel\Sigma$ is locally doubling.
\end{Proposition}

Some of the results of this paper hold for the more general class of rectifiable sets that we now introduce.

\begin{Definition}[Rectifiable sets]\label{def:rectifiable}
We say that a set $R\subset\G$ is  {\em countably $(\G;\G')$-rectifiable} if there exists $\G'$ and  countably many $C^1_H(\G;\G')$-submanifolds $\Sigma_j\subset\G$, 
$j\in\N,$ such that,
denoting by $Q,m$ the homogeneous dimensions of $\G,\G'$, one has 
\[
\psi^{Q-m}\Big(R\setminus\bigcup_{j}\Sigma_j\Big)=0.
\]
We say that $R$ is {\em $(\G;\G')$-rectifiable} if, moreover, $\psi^{Q-m}(R)<+\infty$.
\end{Definition}

The groups $\G,\G'$ will be usually understood and we will simply write {\it rectifiable} in place of {\it $(\G;\G')$-rectifiable}. 
Notice that, by Remark~\ref{rem:V=G'}, if $\psi^{Q-m}(R)>0$, then the group $\G'$ is uniquely determined by $R$ up to biLipschitz isomorphism. 
 We recall also that this notion of rectifiability is not known to be equivalent to the ones by means of cones, as in \cite{FrSeSe2006intrinsic,FranchiSerapioni2016Intrinsic,2019arXiv191200493D,2020arXiv200309196J}.

A key object in the theory of rectifiable sets is the approximate tangent space. 

\begin{Definition}[Approximate tangent space]\label{def:apprtangent}
Let $R\subset\G$ be countably  $(\G;\G')$-rectifiable and let $\Sigma_j$, $j\in\N$, be as in Definition~\ref{def:rectifiable}; for every $\psi^{Q-m}$-a.e.~$p\in R$ we define the {\em approximate tangent space} $T_p^HR$ to $R$ at $p$ as
\[
T^H_p R:=T^H_p\Sigma_{\bar\jmath}\qquad\text{whenever }p\in\Sigma_{\bar\jmath}\setminus\bigcup_{j\le\bar\jmath-1}\Sigma_j.
\]
\end{Definition}

Definition~\ref{def:apprtangent} is well-posed provided one shows that, for $\psi^{Q-m}$-a.e.~$p\in R$, $T^H_p R$ does not change if in Definition~\ref{def:rectifiable} one changes the covering family of submanifolds $(\Sigma_j)_j$. In turn, it is enough to show that,  if $\Sigma',\Sigma''$ are  level sets of  $f'\in C^1_H(\Omega';\G')$,  $f''\in C^1_H(\Omega'';\G')$ defined on  open sets $\Omega',\Omega''\subset\G$ and all points are split-regular  for $f',f''$, then (see also~\cite[Section~2]{DonVittoneJMAA})
\begin{equation}\label{eq:psi of I}
\psi^{Q-m}(\{p\in\Sigma'\cap\Sigma'':T^H_p\Sigma'\neq T^H_p\Sigma''\})=0 .
\end{equation}

Let $I$ be the set in \eqref{eq:psi of I}. Assume by contradiction that $\psi^{Q-m}(I)>0$; we can without loss of generality suppose that $\Sigma'$ is the intrinsic graph of a map $\phi:A\to\V$ defined on an open set $A\subset\W$ for some splitting $\G=\W\cdot\V$. Let $J:=\{w\in A:w\phi(w)\in I\}$; by Theorem~\ref{prop05161504} 
one has $\psi^{Q-m}(J)>0$, hence there exists $\bar w\in J$ such that
\[
\lim_{r\to 0^+}\frac{\psi^{Q-m}(J\cap \Ball(\bar w,r))}{\psi^{Q-m}(\W\cap \Ball(\bar w,r))}=1.
\]
Taking Lemma~\ref{lem09151054} (iii) into account, it is then a routine task to prove that the blow-up of $I$ at $\bar p:=\bar w\phi(\bar w)$, i.e., the limit $\lim_{\lambda\to0^+} \delta_{1/\lambda}(\bar p^{-1}I)$ in the sense of local Hausdorff convergence, is $T^H_{\bar p}\Sigma'$. This implies that $T^H_{\bar p}\Sigma'\supset T^H_{\bar p}\Sigma''$ and in turn, by equality of the dimensions, that $T^H_{\bar p}\Sigma'= T^H_{\bar p}\Sigma''$: this is a contradiction.

\section{The area formula}\label{sec:area}
Let $\PP$ be a homogeneous subgroup of $\G$ with homogeneous dimension $d$ and let $\theta$ be a Haar measure on $\PP$. By dilation invariance of $\scr E$ and $\PP$ one has 
 \begin{align}
	\Theta_{\psi^d}(\theta,0) &= \lim_{\epsilon\to0^+}
		\sup\left\{ \frac{\theta(E\cap\PP)}{\diam(E)^d} : 0\in E\in\scr E,\ 0<\diam(E)\le\epsilon  \right\} \nonumber\\
	&= \lim_{\epsilon\to0^+}
		\sup\left\{ \frac{\theta(\delta_{\diam(E)^{-1}}E\cap\PP)}{\diam(\delta_{\diam(E)^{-1}}E)^d} : 0\in E\in\scr E,\ 0<\diam(E)\le\epsilon \right\}\label{eq2020blabla} \\
	&= \sup\left\{ \theta(E\cap\PP) : 0\in E\in\scr E,\ \diam(E)=1 \right\}.\nonumber
	\end{align}
	This simple observation turns out to be useful to  study the Federer density $\Theta_{\psi^d}$ of $\psi^d\hel \PP$.

\begin{Lemma}\label{lem09122022}
	Let $\PP$ be a homogeneous subgroup of $\G$ with homogeneous dimension~$d$ and let $\psi^d$ be either the spherical or the Hausdorff $d$-dimensional measure on $\G$.
	Then $\psi^d\hel\PP$ is a Haar measure on $\PP$ and for all $x\in \PP$, 
       \begin{equation}\label{eq10041858bis}
		\sup\left\{ \psi^d(E\cap\PP): x\in E\in\scr E,\ \diam(E)=1 \right\}= \Theta_{\psi^d} (\psi^d\hel \PP,x) 
	=  1.
	\end{equation}
\end{Lemma}
\begin{proof}
	As $\scr E$ and $\rho$ are left invariant, $\psi^d\hel \PP$ is a left invariant  measure on $\PP$.
	Therefore, we only need to show that it is non zero and locally finite to prove that it is a Haar measure. Fix a Haar measure $\theta$ on $\PP$. Since $\theta$ is $d$-homogeneous,
	 $\theta$ is Ahlfors $d$-regular on $(\PP,\rho)$, therefore there are constants $0<c<C$ such that 
	\[
	c\theta(B) \le \calH^d(B) \le C\theta(B)
	\]
	for all Borel subsets $B\subset\PP$, see for instance \cite[Exercise~8.11]{Heinonen}.
	By basic comparisons of the Hausdorff and spherical measures, we infer that $\psi^d$ is non zero and locally finite. 
	We can conclude that $\psi^d$ is a Haar measure on $\PP$.
	
	It remains to prove the equalities in~\eqref{eq10041858bis}.
	The first equality now follows from~\eqref{eq2020blabla} and left-invariance.
	The second equality follows instead from Theorem~\ref{thm:Federer}.
\end{proof}

The following lemma proves Theorem~\ref{prop05161504} in a ``linearized'' case and allows to define the area factor $\calA$.

\begin{Lemma}[Definition of the area factor]\label{lem:areafactor}
    Let $\W\cdot\V$ be a splitting of $\G$ with $\W$ normal.
    Assume that $\PP$ is a homogeneous subgroup of $\G$ which is also an intrinsic graph  $\W\to\V$ and let $\Phi_\PP:\W\to\PP$ be the corresponding graph map. Then, there exists a positive constant $\calA(\PP)$, which we call \emph{area factor}, such that 
\[
    \psi^d\hel \PP = \calA(\PP) {\Phi_\PP}_\# (\psi^d\hel \W).
\]
Furthermore, the area factor is continuous in $\PP$.
\end{Lemma}
\begin{proof}
  In order to prove the first part of the lemma it suffices to show that $\mu \defeq {\Phi_\PP}_\# (\psi^d\hel \W)$ is a Haar measure on $\PP$. To see that it is locally finite, note that $\Phi_\PP$ is a homeomorphism between $\W$ and $\PP$ and that therefore bounded open sets in $\PP$ have finite positive  $\mu$ measure. We need to prove that $\mu$ is left invariant. Choose a set $E \subset \PP$.  Let $p = p_\W p_\V$ be a point on $\PP$ and pick a point $x =x_\W x_\V \in E$, we can write
\[
    \pi_\W(px) 
	= \pi_\W(p_\W p_\V x_\W p_\V^{-1} p_\V x_\V )
	= p_\W \varphi(x_\W) ,
\]
where $\varphi:\W\to\W$ is the group automorphism $\varphi(w):= p_\V w p_\V^{-1}$. Let $v\in \g$ be such that $p_\V=\exp(v)$, where $\exp:\g\to\G$ is the exponential map. Then we have
	\[
	\det(\DD\varphi(e)|_\W) = \det(\Ad_{p_\V}|_{\W}) = \det(e^{\ad_v|_{\W}}) 
	= e^{tr(\ad_v|_{\W})} = 1,
	\]
	where $tr(\ad_v|_{\W})=0$ because $\ad_v$ is nilpotent.
	Here, we denoted by $\ad$ and $\Ad$ the adjoint representations of $\g$ and $\G$ respectively; recall that $\Ad_{\exp(v)} = e^{\ad_v}$. This implies that $\varphi$ preserves Haar measures of $\W$ and thus
	\begin{align*}
	\mu(pE)
	= \psi^d (\pi_\W(pE)) 
	= \psi^d(p_\W\varphi(\pi_\W(E))
	= \psi^d(\pi_\W (E)) = \mu(E) .
	\end{align*}
	We conclude that $\mu$ is a Haar measure on $\PP$, so the first part of the statement is proved.

 Let us prove that $\calA(\PP)$ is continuous with respect to $\PP$. By Proposition~\ref{propFedDensity}, $\calA(\PP)$ is equal to $\Theta_{\psi^d}(\mu,0)$ and, by~\eqref{eq2020blabla},
\[
    \calA(\PP)= \sup \{ \psi^d(\pi_\W (E\cap \PP)) : 0\in E\in \scr E, \diam E= 1\}.
\]
Fix $\epsilon>0$ and let $\PP$ and $\PP'$ be homogeneous subgroups that are intrinsic graphs on $\W$ of maps $\phi_\PP,\phi_{\PP'}:\W\to\V$ such that 
\begin{equation}\label{eq12120947}
\rho(\phi_\PP(w),\phi_{\PP'}(w)) <\epsilon
\quad\forall w\in \pi_\W( \cBall(0,1)).
\end{equation}
Pick $E\in \scr E$ with $0\in E$ and $\diam E = 1$ such that $ \psi^d(\pi_\W (E\cap \PP))>(1-\epsilon)\calA(\PP)$. 
Notice that, if $w\in\pi_\W (E\cap \PP)$, then $\rho(w\phi_\PP(w),w\phi_{\PP'}(w)) < \epsilon$.
Therefore, denoting by $\cBall(E,r)$ the closed $r$ neighborhood of $E$, we have
\[
    \pi_\W (E\cap \PP) 
    \subset \pi_\W(\cBall(E,\epsilon)\cap \PP').
\]
If $\psi^d$ is the Hausdorff measure, then $\cBall(E,\epsilon)\in\scr E$ and $\diam(\cBall(E,\epsilon)) \le 1+2\epsilon$;
If $\psi^d$ is the spherical measure, then $E=\cBall(x,1/2)$ for some $x\in\G$ and thus $\cBall(E,\epsilon) \subset \cBall(x,1/2+\epsilon)\in\scr E$ with $\diam(\cBall(x,1/2+\epsilon))\le 1+2\epsilon$.
In both cases, we obtain
\[
     \calA(\PP') \ge (1+2\epsilon)^{-d} (1-\epsilon) \calA( \PP).
\]
Notice that this inequality holds for all $\PP$ and $\PP'$ satisfying \eqref{eq12120947}, hence we also have $\calA(\PP) \ge (1+2\epsilon)^{-d} (1-\epsilon) \calA( \PP')$.
We conclude that $\PP \mapsto \calA(\PP)$ is continuous.
\end{proof}

It is worth observing that the area factor implicitly depends on the fixed group $\W$. We are now ready to prove our first main result.
\begin{proof}[Proof of Theorem~\ref{prop05161504}]
	The function  $a(w) := \calA(T_{w\phi(w)}^H\Sigma) $ is continuous  on $A$ with values in $(0,\infty)$.
	We define the measure $\mu$, supported on $\Sigma$, by
	\[
	\mu(E) := \int_{\pi_\W(E\cap\Sigma)} a(w)\dd (\psi^d\hel \W) (w)
	\]
	for any $E\subset \G$.
	We shall prove~\eqref{eq09071257} by applying Proposition~\ref{propFedDensity},
	that is, we will show that $\Theta_{\psi^d}(\mu;o)=1$ for all $o\in\Sigma$.
	Fix $o\in\Sigma$ and assume without loss of generality that $o=0$.	
	Then
	\begin{align*}
	\Theta_{\psi^d}(\mu;0)
	&= \lim_{r\to 0^+} 
		\sup\left\{ \frac{\mu(E)}{\diam(E)^d} : 0\in E\in\scr E,\ \diam(E)<r\right\} . 
	\end{align*} 
	Using the continuity of the function $a$, we have
	\[
	\Theta_{\psi^d}(\mu;0)
	= a(0) \lim_{r\to0^+} 
		\sup\left\{
		\dfrac{\psi^d(\pi_\W(E \cap \Sigma))}{(\diam E)^d} : 0\in E\in\scr E,\ \diam(E)<r \right\} .
	\]
   Since the projection $\pi_\W$ commutes with dilations, we have for $0<\eta\le 1$,
	\begin{align*}
	\psi^d(\pi_\W(\delta_\eta E\cap\Sigma))= \eta^d \psi^d(\pi_\W( E\cap \delta_{1/\eta} \Sigma)).
	\end{align*}
        Thus 
        \begin{multline*}
	\Theta_{\psi^d}(\mu;0) \\
	= a(0) \lim_{r\to0^+} 
		\sup\left\{ \psi^d(\pi_\W(E \cap \delta_{1/\eta}\Sigma)) : 0\in E\in\scr E,\ \diam(E)=1, 0<\eta <r \right \} .
       \end{multline*}
		We claim that
		\begin{equation}
	\begin{split}\label{eq:tangentdensity}
	\lim_{r \to 0^+} 
		\sup\left\{
		\psi^d(\pi_\W(E\cap\delta_{1/\eta}\Sigma)) : 0\in E\in\scr E,\ \diam(E)=1,\ 0<\eta<r\right\}  \\
	=  \sup\{\psi^d(\pi_\W(E\cap T_0^H\Sigma)) : 0\in E\in \scr E,\ \diam(E)=1 \} .
	\end{split}
	\end{equation}

As in Lemma~\ref{lem09151054}, we denote by $\phi_\eta:\delta_{1/\eta}A\to\V$ the  function whose intrinsic graph is $\delta_{1/\eta} \Sigma$ and by $\phi_0:\W\to\V$ the one for $T_0^H\Sigma$; then, $\phi_\eta$ converges to $\phi_0$ uniformly on compact sets as $\eta\to0$.
In particular, for every $\epsilon>0$ there is $r_\epsilon>0$ such that $\pi_\W(\cBall(0,1))\subset \delta_{1/\eta}A$ and $\rho(\phi_\eta(w),\phi_0(w))<\epsilon$ for all $w\in\pi_\W(\cBall(0,1))$ and $\eta\in(0,r_\epsilon)$.

We start by proving that the left hand side (LHS) of \eqref{eq:tangentdensity} is not greater than the right hand side (RHS); we can assume LHS$>0$.
Fix $\epsilon>0$.
Then there are $\eta\in(0,r_\epsilon)$ and $E$ such that $0\in E\in \scr E$, $\diam E= 1$ and $\psi^d(\pi_\W(E\cap\delta_{1/\eta}\Sigma)) > (1-\epsilon) {\rm LHS}$.
Notice that $\pi_\W(E)\subset\pi_\W(\cBall(0,1))$ and that
\[
      \pi_\W(E\cap\delta_{1/\eta}\Sigma)
      \subset  \pi_\W(\cBall(E,\epsilon)\cap T_0^H\Sigma).
\]
If $\psi^d$ is the Hausdorff measure, then $\tilde E:=\cBall(E,\epsilon)\in\scr E$ and $\diam\tilde E\le 1+2\epsilon$;
If $\psi^d$ is the spherical measure, then $E=\cBall(x,1/2)$ for some $x\in\G$ and thus $\cBall(E,\epsilon)\subset \cBall(x,1/2+\epsilon) =: \tilde E \in \scr E$ and $\diam\tilde E\le 1+2\epsilon$.
Thus, by $d$-homogeneity of $\psi^d \hel \W$, we have
\[
   {\rm RHS} 
   \ge \dfrac{\psi^d(\pi_\W(\tilde E\cap T_0^H\Sigma))}{(\diam \tilde E)^d}
   \ge \dfrac{1-\epsilon}{(1+2\epsilon)^d} {\rm LHS}.
\]
The inequality ${\rm RHS} \ge {\rm LHS}$ follows from the arbitrariness of $\varepsilon$.

For the converse inequality, fix $\epsilon>0$ and $\tilde E$ with $0\in \tilde E\in \scr E$ and $\psi^d(\pi_\W(\tilde E\cap T_0^H\Sigma)) \ge (1-\epsilon){\rm RHS}$. 
Notice that, for every $\eta\in(0,r_\epsilon)$,
\[
	\pi_\W(\delta_{1-2\epsilon}\tilde E\cap T_0^H\Sigma)\subset  \pi_\W(\cBall(\delta_{1-2\epsilon} \tilde E,\epsilon)\cap \delta_{1/\eta}\Sigma)
\]
and that $\diam (\cBall(\delta_{1-2\epsilon} \tilde E,\epsilon))\le 1$. Similarly as before, we can find $\tilde E_\epsilon\in\mathscr E$ such that $\cBall(\delta_{1-2\epsilon} \tilde E,\epsilon)\subset \tilde E_\epsilon$ and $\diam \tilde E_\epsilon= 1$. Therefore, 
\begin{align*}
{\rm LHS} 
		\ge & \limsup_{\eta\to 0^+} 	\psi^d(\pi_\W(\tilde E_\epsilon\cap \delta_{1/\eta}\Sigma))\\
	\ge & \limsup_{\eta\to 0^+} \psi^d(\pi_\W(\cBall(\delta_{1-2\epsilon} \tilde E,\epsilon)\cap \delta_{1/\eta}\Sigma))\\
	\ge &  \psi^d(\pi_\W(\delta_{1-2\epsilon}\tilde E\cap T_0^H\Sigma))\\
	= & (1-2\epsilon)^d\psi^d(\pi_\W(\tilde E\cap T_0^H\Sigma))\\
	\ge &(1-2\epsilon)^d(1-\epsilon){\rm RHS}.
\end{align*}
This concludes the proof of~\eqref{eq:tangentdensity}.

We conclude that, by the definition of the area factor in Lemma~\ref{lem:areafactor}, 
\begin{align*}
\Theta_{\psi^d}(\mu;0)
&= \calA(T_0^H\Sigma) \sup\{\psi^d(\pi_\W(E\cap T_0^H\Sigma)) : 0\in E\in \scr E,\ \diam(E)=1 \} \\
&= \sup\{\psi^d(E\cap T_0^H\Sigma) : 0\in E\in \scr E,\ \diam(E)=1 \}\\
&=1
\end{align*}
where the last equality follows from Lemma~\ref{lem09122022}.
\end{proof}

We conclude this section with some applications of Theorem~\ref{prop05161504}. 
We start by proving the first part in the statement of Corollary \ref{cor:SvsH} about the relation between Hausdorff and spherical Hausdorff measures on rectifiable sets; the second part of Corollary \ref{cor:SvsH}, concerning the same application in the setting of the Heisenberg group endowed with a rotationally invariant distance, will be proved in Proposition~\ref{prop:Hncodim1}

\begin{proof}[Proof of Corollary~\ref{cor:SvsH}, first part]
	If $\PP\in\scr T_{\G,\G'}$, let $\textfrak a(\PP)$ as the constant such that
	\begin{equation}\label{eq03231608}
	\cal S^{Q-m}\hel\PP = \textfrak a(\PP) \cal H^{Q-m}\hel\PP,
	\end{equation}
	which exists because both measures are Haar measures.
	
	Now, let $\G=\W\cdot\V$ be a splitting and $\Sigma$  the $C^1_H$ intrinsic graph of $\phi:A\to\V$ with $A\subset\W$.
		Then, denoting by $\cal A_{\cal S}^\W$ and $\cal A_{\cal H}^\W$ the area factors for the spherical and Hausdorff measures with respect to $\W$,
	\begin{align*}
	\cal S^{Q-m}\hel\Sigma
	&= \cal A_{\cal S}^\W(T^H\Sigma) \Phi_\#(\cal S^{Q-m}\hel\W) \\
	&= \textfrak a(\W)\frac{\cal A_{\cal S}^\W(T^H\Sigma)}{\cal A_{\cal H}^\W(T^H\Sigma)} \cal A_{\cal H}^\W(T^H\Sigma) \Phi_\#(\cal H^{Q-m}\hel\W) \\
	&= \textfrak a(\W)\frac{\cal A_{\cal S}^\W(T^H\Sigma)}{\cal A_{\cal H}^\W(T^H\Sigma)} \cal H^{Q-m}\hel\Sigma .
	\end{align*}
	Since $\Sigma$ is arbitrary, we can apply this equality to $\Sigma=\PP\in\cal W_\W$ to see that
	\[
	\textfrak a(\PP) = \textfrak a(\W)\frac{\cal A_{\cal S}^\W(\PP)}{\cal A_{\cal H}^\W(\PP)} .
	\]
	Continuity of $\textfrak a$ and~\eqref{eq:AAAAA} are now clear.
\end{proof}

\begin{Remark}\label{rem:isodiametric}
The definition of $\textfrak a$ in~\eqref{eq03231608} 
together with Proposition~\ref{propFedDensity} (with $\mu=\cal S^{Q-m}$ and $\psi^d=\cal H^{Q-m}$)  distinctly shows  that the precise value of $\textfrak a(\W)$ is related with the {\it isodiametric problem} on $\W$ about maximizing the measure of subsets of $\W$ with diameter at most 1 (see~\cite{RigotIsodiametric}). This task is a very demanding one already in the Heisenberg group endowed with the Carnot-Carath\'eodory distance, see~\cite{LeonardiRigotV}.
\end{Remark}

We now prove a statement about  weak* convergence of  measures of level sets of $C^1_H$ functions; this will be used in the subsequent Corollary~\ref{cor:densityexistence} as well as later in the proof of the coarea formula.
We note that the proof of Lemma~\ref{lem:weakstarblowup} relies on the Area formula~\eqref{eq09071257}: we are not aware of any alternative strategy.

\begin{Lemma}[Weak* convergence of blow-ups]\label{lem:weakstarblowup}
	Consider an open set $\Omega\subset\G$, a function $g\in C^1_H(\Omega;\G')$ and a point $o\in\Omega$ that is split-regular for $g$. Let $m$ denote the homogeneous dimension of $\G'$ and, for $b\in\G'$ and $\lambda>0$, define 
	\[
	\Sigma_{\lambda,b} :=\delta_{1/\lambda}(o^{-1}\{p\in\Omega:g(p)=g(o)\delta_\lambda b\}) = \{p\in\delta_{1/\lambda}(o^{-1}\Omega):g(o\delta_\lambda p)=g(o)\delta_\lambda b\}.
	\] 
	Then, the weak* convergence of measures
	\[
	\psi^{Q-m}\hel\Sigma_{\lambda,b} 
	\:\overset{*}{\rightharpoonup}\:
	\psi^{Q-m}\hel \{p:\DH  g(o)p=b\} 
	\qquad\text{as }\lambda\to0^+ 
	\]
	holds. 
	Moreover, the convergence is uniform with respect to $b\in\G'$, i.e.,  
	for every $\chi\in C_c(\G)$ and  every $\epsilon>0$ there is $\bar\lambda>0$ such that
	\[
	\left| \int_{\Sigma_{\lambda,b}} \chi \dd\psi^{Q-m} - \int_{\{\DH  g(o)=b\}} \chi \dd\psi^{Q-m}  \right| < \epsilon \qquad\forall\:\lambda\in(0,\bar\lambda), b\in\G'.
	\]
\end{Lemma}
\begin{proof}
	Up to replacing $g$ with the function $x\mapsto g(o)^{-1}g(ox)$, we can assume $o=0$ and $g(o)=0$; in particular, $\Sigma_{\lambda, b}=\delta_{1/\lambda}(\{p\in\Omega:g(p)=\delta_\lambda b\})$.
	Notice that, by Lemma~\ref{lem05161502}, $\Sigma_{\lambda,b}\neq\emptyset$ for all $b$ in a neighborhood of $0$ and  $\lambda$ small enough.
	
	Possibly restricting $\Omega$,
 we can assume that there exists a splitting $\G=\W\cdot\V$, open sets $A\subset\W$, $B\subset\G'$ and a map $\varphi:A\times B\to\V$ such that the statements of Lemma~\ref{lem05161502} hold.
	If $p\in \Sigma_{\lambda,b}$, then there is $a\in A$ such that $p = \delta_{1/\lambda}(a\varphi(a,\delta_\lambda b)) = \delta_{1/\lambda} a \varphi_{\lambda}(\delta_{1/\lambda} a, b)$, 
	where $\varphi_{\lambda}(a,b) := \delta_{1/\lambda}\varphi(\delta_\lambda a,\delta_\lambda b)$.
	In particular, $\Sigma_{\lambda,b}$ is the intrinsic graph of $\varphi_{\lambda}(\cdot,b):\delta_{1/\lambda}A\to\V$.
	
	Denoting by $\varphi_0:\mathbb W\times\G'\to\V$ the implicit function associated with $\DH g(0)$, we have by Lemma~\ref{lem09122015} that
	$\varphi_{\lambda}\to \varphi_0$ uniformly on compact subsets of $\W\times\G'$.
	Moreover
	\begin{align*}
	\lim_{\lambda\to 0^+} T_{a\varphi_{\lambda}(a,b)}\Sigma_{\lambda,b} 
	&= \lim_{\lambda\to 0^+}  T_{\delta_{1/\lambda}(\delta_\lambda a\varphi(\delta_\lambda a,\delta_\lambda b))}\delta_{1/\lambda}\Sigma_{1,\delta_\lambda b} \\
	&= \lim_{\lambda\to 0^+} T_{\delta_\lambda a\varphi(\delta_\lambda a,\delta_\lambda b)} \Sigma_{1,\delta_\lambda b} \\
	&= \lim_{\lambda\to 0^+} \ker(\DH g(\delta_\lambda a\varphi(\delta_\lambda a,\delta_\lambda b))) \\
	&= \ker(\DH  g(0)) 
	\in \scr W_\V,
	\end{align*}
	where the convergence is in the topology of $\scr W_\V$ and it is uniform when $(a,b)$ belong to a compact set of $\W\times\G'$.
	Therefore, using the area formula of Theorem~\ref{prop05161504}, for every $\chi\in C_c(\G)$ we have
	\begin{align}
	\lim_{\lambda\to0^+} \int_{\Sigma_{\lambda,b}} \chi \dd\psi^{Q-m}
	&= \lim_{\lambda\to0^+} \int_{\delta_{1/\lambda}A} \chi(a\varphi_{\lambda}(a,b))\, \calA (T_{a\varphi_{\lambda}(a,b)}\Sigma_{\lambda,b}) \dd \psi^{Q-m}(a) \nonumber\\
	&= \int_{\mathbb W} \chi(a\varphi_0(a,b))\, \calA( \ker\DH  g(0)) \dd \psi^{Q-m}(a) \label{eq:limiteCDD}\\
	&= \int_{\{\DH  g(0)=b\}} \chi \dd\psi^{Q-m}  ,\nonumber
	\end{align}
	where the limit is  uniform when $b$ belongs to a compact subset of $\G'$. Let us show that the convergence is actually uniform on $\G'$.
	
	Since $g$ is Lipschitz continuous in a neighborhood of $0$, there is a positive constant $C$ such that
	$\rho'(0,g(\delta_\lambda p))	\le C\lambda $
	for all $p\in\spt\chi$ and $\lambda$ small enough.
	Therefore, if $\rho'(0,b)>C$, then $\spt\chi\cap \Sigma_{\lambda,b}=\emptyset$.
	Possibly increasing $C$, we can assume that $\spt\chi\cap \{\DH  g(o)=b\} =\emptyset$ for all $b$ such that $\rho'(0,b)>C$. Therefore, the uniformity of the limit~\eqref{eq:limiteCDD} for $b\in\cBall_{\G'}(0,C)$ implies  uniformity  for all $b\in\G'$.
	This completes the proof.
\end{proof}

In the proof of the following corollary, we will need this simple lemma:
\begin{Lemma}\label{lem03231638}
	Let $\theta$ be a Haar measure and $\rho$ a ho\-mo\-ge\-neous  distance on a ho\-mo\-ge\-neous group $\PP$.
	Then $\theta(\partial \Ball_\rho(0,R))=0$ for all $R>0$.
\end{Lemma}
\begin{proof}
By homogeneity, there holds
\begin{align*}
\theta(\partial \Ball(0,R)) 
&= \lim_{\epsilon\to 0^+} \theta( \Ball(0,R+\epsilon)) - \theta( \Ball(0,R-\epsilon))\\ 
&= \theta( \Ball(0,1)) \lim_{\epsilon\to 0^+} ((R+\epsilon)^{\dim_H\PP}- (R-\epsilon)^{\dim_H\PP}) 
=0.
\end{align*}
\end{proof}

\begin{Corollary}\label{cor:densityexistence}
There exists a continuous function $\textfrak d:\scr T_{\G,\G'}\to(0,+\infty)$ with the following property. If $R\subset\G$ is a $(\G;\G')$-rectifiable set and $Q,m$ denote the homogeneous dimensions of $\G,\G'$, respectively, then
\begin{equation}\label{eq:density}
\lim_{r\to 0^+} \frac{\psi^{Q-m}(R\cap\Ball(p,r))}{r^{Q-m}}=\textfrak d(T^H_p R)\qquad\text{for $\psi^{Q-m}$-a.e.~}p\in R.
\end{equation}
Moreover, if $R$ is a $C^1_H$ submanifold, then the equality in~\eqref{eq:density} holds at every $p\in R$.
\end{Corollary}

Clearly, $\textfrak d$  depends  on whether the measure $\psi^{Q-m}$ under consideration is the Hausdorff or the spherical one.

\begin{proof}[Proof of Corollary~\ref{cor:densityexistence}]
Let $\Sigma$ be a $C^1_H$ submanifold and let $\mu:=\psi^{Q-m}\hel(R\setminus\Sigma)$; 
Theorem~\ref{thm:Federer} (ii) implies that
\[
\Theta_{\psi^{Q-m}}(\mu;p)=0\qquad \text{for $\psi^{Q-m}$-a.e.~$p\in\Sigma$},
\]
hence 
\[
\lim_{r\to 0^+} \frac{\psi^{Q-m}((R\setminus\Sigma)\cap\Ball(p,r))}{r^{Q-m}}=0\qquad \text{for $\psi^{Q-m}$-a.e.~$p\in R\cap\Sigma$}.
\]
A similar argument, applied to $\mu:=\psi^{Q-m}\hel(\Sigma\setminus R)$, gives
\[
\lim_{r\to 0^+} \frac{\psi^{Q-m}((\Sigma\setminus R)\cap\Ball(p,r))}{r^{Q-m}}=0\qquad \text{for $\psi^{Q-m}$-a.e.~$p\in R\cap\Sigma$},
\]
i.e.,
\[
\lim_{r\to 0^+} \frac{\psi^{Q-m}(R\cap\Ball(p,r))}{r^{Q-m}}=
\lim_{r\to 0^+} \frac{\psi^{Q-m}(\Sigma\cap\Ball(p,r))}{r^{Q-m}}\qquad \text{for $\psi^{Q-m}$-a.e.~$p\in R\cap\Sigma$}
\]
provided the second limit exists. In particular, it is enough to prove the statement in case $R$ is a $C^1_H$ submanifold.

Let $p\in R$ be fixed; for $\lambda>0$ define $R_\lambda:=\delta_{1/\lambda}(p^{-1}R)$ and, by Lemma~\ref{lem:weakstarblowup}, $\psi^{Q-m}\hel R_\lambda \:\overset{*}{\rightharpoonup}\: 	\psi^{Q-m}\hel T^H_p R $. Since $\psi^{Q-m}(T^H_pR \cap\partial\Ball(0,1))=0$, 
using~\cite[Proposition 1.62 (b)]{AmFuPa2000} and Lemma~\ref{lem03231638}, one gets
\[
\lim_{r\to 0^+} \frac{\psi^{Q-m}(R\cap\Ball(p,r))}{r^{Q-m}}
= \lim_{r\to 0^+} \psi^{Q-m}(R_\lambda\cap\Ball(0,1))
= \psi^{Q-m}(T^H_pR \cap\Ball(0,1)).
\]
Statement~\eqref{eq:density} follows on setting $\textfrak d(\PP):=\psi^{Q-m}(\PP \cap\Ball(0,1))$ for every $\PP\in\scr T_{\G,\G'}$.

It remains only to prove
the continuity of $\textfrak d$ at every fixed $\W\in\scr T_{\G,\G'}$. Every $\PP\in\scr T_{\G,\G'}$ in a proper neighborhood of $\W$ is an intrinsic graph over $\W$. Denoting by $\pi_{\W}:\G\to\W$ the projection defined in \eqref{eq:defproiezioni}, we have by Lemma~\ref{lem:areafactor} that
\[
\textfrak d(\PP)=\psi^{Q-m}(\PP \cap\Ball(0,1)) = \calA(\PP)\psi^{Q-m}(\pi_\W(\PP \cap\Ball(0,1))),
\]
hence we have to prove only the continuity of $\PP\mapsto \psi^{Q-m}(\pi_\W(\PP \cap\Ball(0,1)))$ at $\W$. Let $\epsilon>0$ be fixed; then, if $\PP$ is close enough to $\W$, one has
\[
\W\cap\Ball(0,1-\epsilon)\subset \pi_\W(\PP \cap\Ball(0,1))\subset \W\cap\Ball(0,1+\epsilon)
\]
and the continuity of $\PP\mapsto \psi^{Q-m}(\pi_\W(\PP \cap\Ball(0,1)))$ at $\W$ follows.
\end{proof}

We conclude this section with the following result, similar in spirit to Lemma~\ref{lem:weakstarblowup}.
It will be used in the proof of Lemma~\ref{lem:coareacont}.

\begin{Corollary}\label{corollaryreplacement}
	Suppose that, for $n\in \N$, $L_n:\G\to \G'$ is a homogeneous morphism and that the $L_n$ converge to a surjective homogeneous morphism $L:\G\to \G'$ such that $\ker L$ splits $\G$. 
	Then the following weak* convergence of measures holds: 
	  \[
	    	\psi^{Q-m}\hel \{L_n= s\}
		\overset{*}{\rightharpoonup}
		\psi^{Q-m}\hel \{L=s\} ,
	  \]
	where $Q$ is the homogeneous dimension of $\G$ and $m$ is the homogeneous dimension of $\G'$. 
	More precisely, given a function $\chi \in C_c(\G)$ and $\epsilon>0$, 
	there exists $N\in\N$ such that for all $n\ge N$ and $s\in \G'$
	\[
	    	\left \vert \int_{\{L_n= s\}} \chi \dd  \psi^{Q-m} - 
		\int_{ \{L= s\}} \chi \dd \psi^{Q-m}\right \vert
	      <\epsilon.
	\]
\end{Corollary}
\begin{proof}
	Denote by $\W:=\ker L$ and let $\G=\W\cdot\V$ a splitting.	
	Recall that $\V$ and $\G'$ are also vector spaces, the morphisms $L_n$ are linear maps 
	and that $L|_\V:\V\to\G'$ is an isomorphism.
	Therefore, there exists $N\in\N$ such that $L_n|_{\V}$ is  an isomorphism for all $n\ge N$.
	For all such $n$ and $s\in\G'$, define $\phi_n^s:\W\to\V$ by
	\[
	\phi_n^s(w) := L_n|_{\V}^{-1}(L_n(w)^{-1}s) .
	\]
	Notice that $\{L_n=s\}$ is the intrinsic graph of $\phi_n^s$.
	Let $\phi_\infty^s:\W\to\V$ the function whose intrinsic graph is $\{L=s\}$:
	it is clear that $\phi_n^s(w)\to \phi_\infty^s(w)$ uniformly on compact sets in the variables $(w,s)\in\W\times\G'$.

	Fix $\chi\in C_c(\G)$.
	Then
	\[
	\int_{\{L_n=s\}} \chi \dd \psi^{Q-m}
	= \cal A(\ker L_n)  \int_\W \chi(w\phi_n^s(w)) \dd\psi^{Q-m}(w) 
	\]
	where the functions $\tilde\chi_n:(s,w)\mapsto \chi(w\phi_n^s(w))$ are continuous and uniformly converge to $ (s,w)\mapsto(wL|_{\V}^{-1}(s))$ as $n\to\infty$.
	Moreover, $\cal A(\ker L_n)\to 1$.
	This completes the proof.
\end{proof}

\section{The coarea formula}
\subsection{Set-up}

Let $\G$ be a Carnot group, $\rho$ a homogeneous distance on $\G$ and $Q$ the homogeneous dimension of $\G$. 
Let also $\MM$, $\LL$ and $\KK$ be graded groups and 
 such that  $\LL\MM = \KK$ and $\MM\cap\LL=\{0\}$;
let $m$ and $\ell$ be the respective homogeneous dimensions of $\MM$ and $\LL$.

Our aim is to prove Theorem~\ref{thm:coarea}, which by Proposition~\ref{propFedDensity} will be a consequence of the following Theorem~\ref{thm:coarea2}: here, $\calC(\PP, L)$  denotes the \emph{coarea factor} corresponding to a homogeneous subgroup $\PP$ of $\G$ and a homogeneous morphism $L:\G\to \LL$; the coarea factor is going to be defined later in Proposition \ref{prop:coarea-factor}. The function $\calC(\PP,L)$ is continuous in $\PP$ and $L$, see Lemma~\ref{lem:coareacont}.

\begin{Theorem}\label{thm:coarea2}
	Let $\Omega\subset\G$ be open, let $f\in C^1_H(\Omega;\MM)$ and  assume that all points in $\Omega$ are split-regular for $f$, so that  $\Sigma:=\{p\in\Omega:f(p)=0\}$ is a $C^1_H$ submanifold. Consider a function $u:\Omega\to\LL$  such that $uf\in C^1_H(\Omega;\KK)$ and assume that
\begin{equation*}
\text{for $\psi^{Q-m}$-a.e.~$p\in\Sigma$,}\quad
\left\{
\begin{array}{l}
\text{either $\DH  (uf)_p|_{T^H_p\Sigma}$ is not surjective on $\LL$,}\\
\text{or $p$ is split-regular for $uf$.}
\end{array}
\right.
\end{equation*}
	For $s\in\LL$ set $\Sigma^s:=\Sigma\cap u^{-1}(s)$. Then
	\begin{itemize}
	\item[(i)] for every Borel set $E\subset \Omega$ the function $\LL\ni s\mapsto\psi^{Q-m-\ell}(E\cap \Sigma^s)\in[0,+\infty]$ is $\psi^\ell$-measurable;
	\item[(ii)] the  function
	\begin{equation}\label{eq:defmuSigmau}
             \mu_{\Sigma,u}(E) \defeq \int_{\LL}  \psi^{Q-m-\ell}(E\cap \Sigma^s) \dd \cal\psi^\ell (s),
        \end{equation}
defined on Borel sets, is a   locally finite measure; 
	\item[(iii)] the Radon--Nikodym density $\Theta$ of $\mu_{\Sigma,u}$ with respect to $\psi^{Q-m}\hel\Sigma$ of  is locally bounded and
	\begin{equation}\label{eq:coarealocal}
             \Theta(p)= \calC (T^H_p\Sigma,\DH  (uf)_p)\qquad\text{for $\psi^{Q-m}$-a.e.~}p\in\Sigma.
	\end{equation}
\end{itemize}	
\end{Theorem}

\begin{Remark}\label{rem:olyonu}
Let us prove that the differential $\DH  (uf)_p|_{T^H_p\Sigma}$ depends only on the restriction of $u$ to $\Sigma$ and, moreover, that it does not depend on the choice of the defining function $f$ for $\Sigma$.
In particular, in view of Proposition~\ref{prop:coarea-factor} also the coarea factor $\calC (T^H_p\Sigma,\DH  (uf)_p)$ depends only on the restriction of $u$ to $\Sigma$. 

Let $v\in T^H_p\Sigma$; then there exist sequences $r_j\to 0^+$ and $q_j\to p$ such that $q_j\in\Sigma$ and $v=\lim_{j\to\infty}\delta_{1/r_j} (p^{-1}q_j)$. In particular, $\|q_j^{-1}p\delta_{r_j}v\|_\rho=o(r_j)$ and, by Lemma~\ref{lem10051109}, 
\[
\lim_{j\to\infty}\delta_{1/r_j}\left( (uf)(q_j)^{-1}(uf)(p\delta_{r_j}v)\right)=0.
\]
Since $f|_\Sigma=0$ we obtain
\begin{align*}
\DH (uf)_p(v) &= \lim_{j\to\infty}\delta_{1/r_j}\left( (uf)(p)^{-1}(uf)(p\delta_{r_j}v)\right)\\
&= \lim_{j\to\infty}\delta_{1/r_j}\left( (uf)(p)^{-1}(uf)(q_j)\right)\\
&= \lim_{j\to\infty}\delta_{1/r_j}\left( u(p)^{-1}u(q_j)\right).
\end{align*}
This proves what claimed.
\end{Remark}

The proof of Theorem~\ref{thm:coarea2} is divided into several steps. We start by proving that $\mu_{\Sigma,u}$ is a well defined locally finite measure concentrated on $\Sigma$; this uses an abstract coarea inequality. Then we consider the linear case in order to  apply a blow-up argument; in doing so, we will define the coarea factor. We finally consider separately  ``good points'', i.e., those where $\DH (uf)|_{T^H\Sigma}$  has full rank, and  ``bad points'', where $\dD_H (uf)\vert_{T^H \Sigma}$ is not surjective:  at good points the blow-up argument applies, while  the set of bad points is negligible  by an argument similar to the proof of the coarea inequality.

\subsection{Coarea Inequality}
In this section we prove  Proposition~\ref{prop07171752}, which is a  consequence of the following Lemma~\ref{lem-coarea-ineq}; the latter  is basically \cite[Theorem~2.10.25]{FedererGMT}, with a slightly different use of the Lipschitz constant. See also~\cite[Theorem 1.4]{MagnaniCoareaInequality} and~\cite[Lemma~3.5]{EvansGariepyRevised}.

\begin{Lemma}[Abstract Coarea Inequality]\label{lem-coarea-ineq}
	Let $(X,d_X)$ and $(Y,d_Y)$ be boundedly compact metric spaces and assume that there exist $\beta\ge0$ and $C\ge 0$ such that
	\[
	\cal H^\beta(E) \le C \diam(E)^\beta\qquad \text{for all } E\subset Y,
	\]
	where $\cal H^\beta$ is the $\beta$-dimensional Hausdorff measure on $(Y,d_Y)$.
	Let $u:X\to Y$ be a locally Lipschitz function and for $\epsilon>0$ consider
	\[
	\Lip_\epsilon(u) := \sup\left\{\frac{d_Y(u(x),u(y))}{d_X(x,y)} : 0<d_X(x,y)<\epsilon \right\},\qquad \Lip_0(u) := \lim_{\epsilon\to0} \Lip_\epsilon(u).
	\]
	Then, for every $\alpha\ge\beta$ and  every Borel set $A\subset X$ with $\cal H^\alpha(A)<+\infty$,
	the function $y\mapsto \cal H^{\alpha-\beta}(u^{-1}(y)\cap A)$ is $\cal H^\beta$-measurable and
	\begin{equation}\nonumber
	\int_Y \cal H^{\alpha-\beta}(u^{-1}(y)\cap A) \dd \cal H^\beta(y)
	\le C \Lip_0(u)^\beta \cal H^\alpha(A) .
	\end{equation}
	Moreover, the set function $A\mapsto \int_Y \cal H^{\alpha-\beta}(u^{-1}(y)\cap A) \dd \cal H^\beta(y)$ is a Borel measure.
\end{Lemma}

The proof is standard. 
In our setting, the ``abstract'' coarea inequality translates as follows.

\begin{Proposition}[Coarea inequality]\label{prop07171752}
	Under the assumptions and notation of  Theorem~\ref{thm:coarea2}, one has 
\begin{itemize}
\item[(i)] $u|_\Sigma$ is locally Lipschitz continuous;
\item[(ii)] for every Borel set $E\subset\G$, the function $\LL\ni s\mapsto\psi^{Q-m-\ell}(E\cap \Sigma^s)\in[0,+\infty]$ is $\psi^\ell$-measurable;
\item[(iii)] for every compact $K\subset\Sigma$, the  coarea inequality 
	\[
	\mu_{\Sigma,u}(K)
	\le C \Lip(u|_K)^\ell \psi^{Q-m}(K)
	\]
	holds for a suitable $C=C(\LL)>0$;
\item[(iv)] $\mu_{\Sigma,u}$ is a Borel measure on $\Omega$, $\mu_{\Sigma,u}\ll\psi^{Q-m}\hel\Sigma$ with locally bounded density.
\end{itemize}
	
\end{Proposition}
\begin{proof}
The local Lipschitz continuity of  $u|_{\Sigma}$ follows from Lemma~\ref{lem10051109} because of the assumption $uf\in C^1_H(\Omega;\KK)$ and the fact that $u|_\Sigma=uf|_\Sigma$. As already noticed in the proof of Lemma~\ref{lem-coarea-ineq}, statement (ii) follows from~\cite[2.10.26]{FedererGMT}; the careful reader will observe that~\cite[2.10.26]{FedererGMT} is stated only when $\psi^\ell=\cal H^\ell$, but its proof easily adapts to the case $\psi^\ell=\cal S^\ell$. Concerning statement (iii), we notice that $\psi^{Q-m}(K)<\infty$ because the measure $\psi^{Q-m}\hel\Sigma$ is locally finite by Lemma~\ref{lem09151054} and the area formula (Theorem~\ref{prop05161504}): in particular, one can apply Lemma~\ref{lem-coarea-ineq}. Statement (iv) is now a consequence of statement (iii) and the Radon--Nikodym Theorem, which can be applied because $\psi^{Q-m}\hel\Sigma$ is doubling by~\eqref{eq:Ahlfors}  (see, e.g.,~\cite{rigot2018differentiation}).
\end{proof}

\subsection{Linear case: definition of the coarea factor}
In the following Proposition~\ref{prop:coarea-factor} we prove the coarea formula in a ``linear'' case, and in doing so we will introduce the coarea factor. We are going to consider a homogeneous subgroup $\PP$ of $\G$ that is also a $C^1_H$ submanifold. We observe that this implies that $\PP$ coincides with its homogeneous tangent subgroup; in particular, $\PP$ is normal and it is the kernel of a surjective homogeneous morphism on $\G$.

\begin{Proposition}[Definition of coarea factor]\label{prop:coarea-factor}
	Let $\PP$ be a homogeneous  subgroup of $\G$ that is a $C^1_H$ submanifold of $\G$ and let $L:\PP\to\LL$ be a homogeneous morphism.  Let $\mu_{\PP,L}$ be as in~\eqref{eq:defmuSigmau}, namely,
\[
\mu_{\PP,L}  
	\defeq \int_{\LL} \psi^{Q-m-\ell}\hel  L^{-1}(s) \dd \psi^\ell(s).
\]	
	Then, $\mu_{\PP,L}$ is either null or a Haar measure on $\PP$. In particular, there  exists  $\cal C (\PP,L)\ge0$, which we call \emph{coarea factor}, such that	
	\begin{equation}\label{eq07171719}
	\mu_{\PP,L}  
	= \calC(\PP,L)\ \psi^{Q-m} \hel \PP.
	\end{equation}
	Moreover, $\cal C(\PP,L)>0$ if and only if $L(\PP)=\LL$.
	\end{Proposition}
\begin{proof}
	Since $L$ is Lipschitz on $\PP$, 
	we can apply Proposition~\ref{prop07171752}
	and obtain that $\mu_{\PP,L}$
	is a well defined Borel regular measure that is also 
	absolutely continuous with respect to $\psi^{Q-m}\hel \PP$ and finite on bounded sets.
	
	If $L(\PP)\neq \LL$,
	then $\mu_{\PP,L}=0$ and thus~\eqref{eq07171719} holds with $\cal C(\PP,L)=0$.
	
	If $L(\PP)=\LL$, then we will show that $\mu_{\PP,L}$ is a Haar measure on $\PP$, which is equivalent to \eqref{eq07171719} with
        $\cal C(\PP,L)> 0$.
        For $s\in\LL$ let $\PP^s:=L^{-1}(s)$.
		Since $\PP^s$ is a coset of $\PP^0$, 
	$\psi^{Q-m-\ell}\hel \PP^s$ 
	is the push-forward of $\psi^{Q-m-\ell}\hel \PP^0$ (which is a Haar measure on $\PP^0$) via a left translation.
	It follows that $\mu_{\PP,L}$ is nonzero on nonempty open subsets of $\PP$.
	
	\noindent We need only to show that $\mu_{\PP,L}$ is left-invariant:
	let $p\in \PP$ and choose a Borel set $A\subset \PP$.	
	For every $s\in\LL$ we have
	$p^{-1} \PP^s = \{q\in \PP: L(pq) = s\} = \PP^{L(p)^{-1}s}$.
	By left invariance of $\psi^{Q-m-\ell}$ and $\psi^\ell$, we have
	\begin{align*}
	    \mu_{\PP,L}(p A) &= \int_{\LL} \psi^{Q-m-\ell}((pA)\cap \PP^s)\dd \psi^{\ell}(s) \\
	    &= \int_{\R^k} \psi^{Q-m-\ell}(p(A\cap \PP^{L(p)^{-1}s}))\dd \psi^{\ell}(s) \\
	    &= \int_{\R^k} \psi^{Q-m-\ell}(A\cap \PP^{L(p)^{-1}s})\dd \psi^{\ell}(s) \\
	    &= \mu_{\PP,L}(A)
	\end{align*}
	as wished.
\end{proof}

We now prove a continuity property for the coarea factor $\cal C(\PP,L)$. We agree that, when $L:\G\to\LL$ is defined on the whole $\G$, the symbol $\cal C(\PP,L)$ stands for $\cal C(\PP,L|_\PP)$.

\begin{Lemma}\label{lem:coareacont}
Assume that, for $n\in\N$,   surjective homogeneous morphisms $F,F_n:\G\to\MM$ and homogeneous maps $L,L_n:\G\to\LL$ are given in such a way that
\begin{enumerate}[label=(\roman*)]
\item $LF$ and $L_nF_n$ are homogeneous morphisms $\G\to\KK$;
\item $\ker F$ and $\ker(LF)$ split $\G$;
\item $F_n\to F$ and $L_n\to L$  on $\G$ as $n\to\infty$.
\end{enumerate}
Then $\cal C(\ker F_n,L_n)\to\cal C(\ker F,L)$ as $n\to\infty$.
\end{Lemma}

\begin{proof}
Set $\PP_n:=\ker F_n$ and $\PP:=\ker F$; let $\V$ be a complementary subgroup to $\PP$. Then, for large enough $n$, $\PP_n\cdot\V$ is a splitting of $\G$ and  the subgroup $\PP_n$ is the intrinsic graph $\PP\to\V$ of a homogeneous map $\phi_n\in C^1_{\PP,\V}(\PP)$. Observe that $\phi_n\to 0$ locally uniformly on $\PP$ because $\PP_n\to\PP$.
This, together with Lemma~\ref{lem:areafactor}  and the continuity of the area factor by Lemma~\ref{lem:areafactor}, implies that $\psi^{Q-m}\hel \PP_n$ converges weakly* to $\psi^{Q-m}\hel \PP$. Therefore, by Proposition~\ref{prop:coarea-factor} we have only to show that
\begin{equation}\label{eq:weakconvint}
    \mu_{\PP_n,L_n} \stackrel{*}{\rightharpoonup} \mu_{\PP,L}.
\end{equation}
   
If $L|_\PP$ is surjective, also $LF$ is surjective.
Since $\ker(LF)$ splits $\G$, then~\eqref{eq:weakconvint} follows from Corollary~\ref{corollaryreplacement}.
If $L\vert_\PP$ is not surjective, we can without loss of generality suppose that $L_n\vert_{\PP_n}$ is surjective for all $n$. By homogeneity, it suffices to prove that $\mu_{\PP_n,L_n}(\cBall_\G(0,1)) \to 0$. We have
\begin{align*}
      \mu_{\PP_n,L_n}(\cBall_\G(0,1))&=\int_{\LL} \psi^{Q-m-\ell}(\PP_n\cap L_n^{-1}(s)\cap\cBall_\G(0,1))\dd \psi^{\ell}(s)\\
&\le \psi^{\ell}(L_n(\cBall_\G(0,1)\cap \PP_n)) \sup_{s\in\LL} \psi^{Q-m-\ell}(\cBall_\G(0,1)\cap \PP_n\cap L_n^{-1}(s))\\
&\le \psi^{\ell}(L_n(\cBall_\G(0,1)\cap \PP_n))
\end{align*}
where the last inequality holds because, considering $p\in\PP_n$  such that $L_n(p)=s^{-1}$, we have by  Lemma~\ref{lem09122022}
\[
\psi^{Q-m-\ell}(\cBall_\G(0,1)\cap \PP_n\cap L_n^{-1}(s))=\psi^{Q-m-\ell}(\cBall_\G(p,1)\cap \PP_n\cap L_n^{-1}(0))\le 2^{Q-m-\ell}.
\]
Thus, we have to prove that $\psi^{\ell}(L_n(\cBall_\G(0,1)\cap \PP_n))\to0$; notice that $L_n(\cBall_\G(0,1)\cap \PP_n)$ converges in the Hausdorff distance to $L(\cBall_\G(0,1)\cap \PP)$, which is a compact set contained in a strict subspace of $\LL$. As $\psi^\ell$ is a Haar measure on $\LL$, we have $\psi^{\ell}(L_n(\cBall_\G(0,1)\cap\PP_n)\to 0$ as $n\to \infty$. 
  \end{proof}

\subsection{Good Points}\label{sec:goodpoints}
By ``good'' point $o\in\Sigma$ we mean a point where the differential $\DH (uf)|_{T_o^H\Sigma}$ is surjective onto $\LL$; the following Proposition~\ref{prop07251509} shows that the Radon--Nikodym density $\Theta$ of $\mu_{\Sigma,u}$ with respect to $\psi^{Q-m}\hel\Sigma$ can be explicitly computed at its Lebesgue points and coincides with the coarea factor. 
Notice that almost every $o\in\Sigma$ is a Lebesgue point for $\Theta$, in the sense that
\begin{equation}\label{eq:RNo}
\lim_{r\to 0^+}\aveint{\Sigma\cap\Ball(o,r)}{} \left|\Theta - \Theta(o)\right|\dd \psi^{Q-m}=0. 
\end{equation}
 
\begin{Proposition}\label{prop07251509}
	Under the assumptions and notation of Theorem~\ref{thm:coarea2}, one has that the equality
	\begin{equation}\label{eq:Theta=coarea}
    \Theta(o) = \calC(T^H_o\Sigma,\DH (uf)(o)).
    \end{equation}
	holds for $\psi^{Q-m}$-a.e.~$o\in\Sigma$ such that  $\DH (uf)\vert_{T^H_o \Sigma}$ is onto $\LL$.
\end{Proposition}
\begin{proof}
	We are going to prove~\eqref{eq:Theta=coarea} for all $o\in\Sigma$ such that  $\DH (uf)\vert_{T^H_0 \Sigma}$ is onto $\LL$, $o$ is split-regular for $uf$ and~\eqref{eq:RNo} holds; up to left translations, we may assume that $o=0$ and $u(0)=0$. 
	For every Borel set $A\subset\G$ and $\lambda>0$ we have, on the one hand
	\begin{align*}
	\mu_{\Sigma,u}(\delta_\lambda A) 
	&= \int_{\Sigma\cap\delta_\lambda A} \Theta(p) \dd \psi^{Q-m}(p) \\
	&= \lambda^{Q-m} \int_{(\delta_{1/\lambda}\Sigma)\cap A} \Theta(\delta_\lambda p) \dd \psi^{Q-m}(p) \\
	&= \lambda^{Q-m} (\Theta\circ\delta_\lambda) \psi^{Q-m}\hel\delta_{1/\lambda}\Sigma (A) .
	\end{align*}
	On the other hand,
	\begin{align*}
	\mu_{\Sigma,u}(\delta_\lambda A) 
	&= \int_{\LL} \psi^{Q-m-\ell}((\delta_\lambda A)\cap \Sigma \cap \{u=s\}) \dd \psi^\ell(s) \\
	&= \int_{\LL} \psi^{Q-m-\ell}(\delta_\lambda (A\cap \delta_{1/\lambda}\Sigma \cap \{u_\lambda=\delta_{1/\lambda}s\})) \dd \psi^\ell(s) \\
	&= \lambda^{Q-m} \int_{\LL} \psi^{Q-m-\ell}(A\cap \delta_{1/\lambda}\Sigma \cap \{u_\lambda=t\})) \dd \psi^\ell(t) , \\ 
	\end{align*}
	where $u_\lambda(p) := \delta_{1/\lambda}u(\delta_\lambda p)$.
	Therefore, one has the equality of measures
	\begin{equation}\label{eq09181448}
	(\Theta\circ\delta_\lambda) \psi^{Q-m}\hel\delta_{1/\lambda}\Sigma
	=  \int_{\LL} \psi^{Q-m-\ell}\hel (\delta_{1/\lambda}\Sigma \cap \{u_\lambda=b\}) \dd \psi^\ell(b) .
	\end{equation}
	We now compute the weak* limits as $\lambda\to 0^+$ of  each side of~\eqref{eq09181448}. 
	Concerning the left-hand side, for every $\chi\in C_c(\G)$ one has
	\begin{multline*}
	\int_{\delta_{1/\lambda}\Sigma } \chi(p) \Theta(\delta_\lambda p)  \dd\psi^{Q-m}(p)\\
	= \int_{\delta_{1/\lambda}\Sigma } \chi(p) (\Theta(\delta_\lambda p)-\Theta(0))  \dd\psi^{Q-m}(p)
		+ \Theta(0) \int_{\delta_{1/\lambda}\Sigma } \chi(p) \dd\psi^{Q-m}(p),
	\end{multline*}
Let $r>0$ be such that $\spt\chi\subset\Ball(0,r)$, then
	\begin{align*}
	&\hspace{-2cm}\left| \int_{\delta_{1/\lambda}\Sigma } \chi(p) (\Theta(\delta_\lambda p)-\Theta(0))  \dd\psi^{Q-m} (p)\right| \\
	\le& \|\chi\|_{\infty}  \int_{\Ball(0,r)\cap\delta_{1/\lambda}\Sigma } |\Theta(\delta_\lambda p)-\Theta(0)|  \dd\psi^{Q-m}(p)  \\
	=& \|\chi\|_{\infty} \lambda^{m-Q} \int_{\Ball(0,\lambda r)\cap\Sigma} |\Theta(p)-\Theta(0)|  \dd\psi^{Q-m}(p) \\
	\le&  C\:\|\chi\|_{\infty} \aveint{\Ball(0,\lambda r)\cap\Sigma}{} |\Theta(p)-\Theta(0)|  \dd\psi^{Q-m}(p)
	\end{align*}
	for a suitable positive $C$. Exploiting~\eqref{eq:RNo} one gets
	\begin{equation}\label{eq:LHS}
	\begin{split}
	\lim_{\lambda\to0^+}\int_{\delta_{1/\lambda}\Sigma } \chi(p)\Theta(\delta_\lambda p) \dd\psi^{Q-m}(p)
	&=\Theta(0)\lim_{\lambda\to0^+}\int_{\delta_{1/\lambda}\Sigma } \chi(p)   \dd\psi^{Q-m}(p)\\
	&= \Theta(0)\int_{T^H_0\Sigma} \chi(p) \dd\psi^{Q-m}(p),
	\end{split}
	\end{equation}
	the last equality following from Lemma~\ref{lem:weakstarblowup}.
	
	We now consider the right-hand side of~\eqref{eq09181448}; setting $(uf)_\lambda(p):=\delta_{1/\lambda}((uf)(\delta_\lambda p))$,  for every $\chi\in C_c(\G)$ one has
	\begin{equation}\nonumber
	\begin{split}
	 \lim_{\lambda\to 0^+}\int_\LL \int_{\delta_{1/\lambda}\Sigma\cap \{u_\lambda=b\}}\chi\dd\psi^{Q-m-\ell}\dd\psi^\ell(b)
= & \lim_{\lambda\to 0^+}\int_\LL \int_{\{(fu)_\lambda=b\}}\chi\dd\psi^{Q-m-\ell}\dd\psi^\ell(b)\\
= & \int_\LL \int_{ \{\DH (uf)(0)=b\} }\chi\dd\psi^{Q-m-\ell}\dd\psi^\ell(b),
	\end{split}
	\end{equation}
where we used Lemma~\ref{lem:weakstarblowup}. The definition of coarea factor then gives
	\begin{equation}\label{eq:RHS}
	\begin{split}
\lim_{\lambda\to 0^+}\int_\LL \int_{\delta_{1/\lambda}\Sigma\cap \{u_\lambda=b\}}\chi\dd\psi^{Q-m-\ell}\dd\psi^\ell(b)=
 & \int \chi\dd\mu_{T^H_0\Sigma,\DH (uf)(0)}\\
 &\hspace{-1cm}= \calC(T^H_0\Sigma,\DH (uf)(0))\int_{T^H_0\Sigma}\chi\dd\psi^{Q-m}
	\end{split}
	\end{equation}
The statement is now a consequence of~\eqref{eq09181448},~\eqref{eq:LHS} and~\eqref{eq:RHS}.
\end{proof}

\subsection{Bad points}\label{sec:badpoints}
In contrast with ``good'' ones, ``bad'' points are those points $p$ where $(\DH (uf))(p)|_{T^H_p\Sigma}$ is not surjective. The following lemma states that they are $\mu_{\Sigma,u}$-negligible: a posteriori, this is consistent with the fact that, by definition,  the coarea factor is null at such points.

\begin{Lemma}\label{lem07241228}
	Under the assumptions and notation of Theorem~\ref{thm:coarea2}, one has
	\[
	\mu_{\Sigma,u}(\{p\in\Sigma:\DH (uf)(p)|_{T^H_p\Sigma}\text{ is not onto }\LL\})=0.
	\]
\end{Lemma}
\begin{proof}
It is enough to show that $\mu_{\Sigma,u}(E)=0$ for an arbitrary compact subset $E$ of $\{p\in\Sigma:\DH (uf)(p)|_{T^H_p\Sigma}\text{ is not onto }\LL\}$, which is closed. 
We have $\psi^{Q-m}(E)<\infty$. 
Fix $\epsilon>0$;
by the compactness of $E$ and the locally uniform differentiability of both $f$ and $uf$, 
there exists $r>0$ such that 
$\cBall(E,r)\subset\Omega$ and,
for all $p\in E$ and all $q\in \Sigma\cap \Ball(p,r)$, the inequalities
\[
      \dist(q,p T^H_p\Sigma) \le \epsilon \rho_\G(p,q) ,
\]
and
\[
	\rho_\KK \left( \DH (uf)_p(p^{-1}q) , (uf)(p)^{-1}(uf)(q) \right) 
	\le M\epsilon \rho_\G(p , q )
\]
hold,
where $M=\Lip((uf)|_{\cBall(E,r)})$.
Fixing a positive integer $j>1/r$,
one can cover $E$ by countably many closed sets $\{B_i^j\}_i$ of diameter $d^j_i:=\diam B^j_i$ belonging to the class $\scr E$ and  such that
\begin{equation}\label{eq:quasiHausdorff}
d^j_i<1/j,\text{ for all $i$,}\quad\text{ and}\quad \sum_i (d_i^j)^{Q-m} <\psi^{Q-m}(E)+ 1/j.
\end{equation}
Imitating the proof of \cite[Lemma~3.5]{EvansGariepyRevised},
we define the functions $g^j_i:\LL\to[0,1]$ by $g^j_i = (d^j_i)^{Q-m-\ell} \mathbf{1}_{u(B_i^j\cap\Sigma)}$. 
	Note that, using the standard notation $\psi^{Q-m-\ell}_\delta$ for the  pre-measures used in the Carathéodory construction, one has
	\begin{equation}\label{eq07171643}
	\psi^{Q-m-\ell}_{1/j}(u^{-1}(y)\cap E) \le \sum_i g^j_i(y) ,
	\end{equation}
	for all $y\in Y$. 

Then one gets, using upper integrals,
\begin{equation}\label{eq03261658}
\begin{aligned}
	\int_{\LL} \psi^{Q-m-\ell}_{1/j}(E\cap u^{-1}(s)) \dd \psi^\ell(s)
		&\overset{\eqref{eq07171643}}{\le} \int_\LL \sum_i g^j_i(y) \dd \psi^\ell(y) \\
	&\overset{*}{\le} \sum_i \int_\LL  g^j_i(y) \dd \psi^\ell(y) \\
	&\le  \int_\LL\sum_i(d^j_i)^{Q-m-\ell}\mathbf{1}_{u(B^j_i\cap\Sigma)}(s)\dd\psi^\ell(s)\\
	& \le \sum_i (d_i^j)^{Q-m-\ell} \psi^{\ell}(u(B_j^i\cap \Sigma)) ,
\end{aligned}
\end{equation}
where the inequality marked by $*$ follow from Fatou's Lemma.
We claim that 
\begin{equation}\label{eq:volumecontrol}
  \psi^{\ell}(u(B_i^j\cap\Sigma))\le M^\ell C(\epsilon,\LL) (\diam B_i^j)^\ell,
\end{equation}
for a suitable $C(\epsilon,\LL)>0$ such that $\lim_{\epsilon\to 0^+}C(\epsilon,\LL)=0$. 

Let us prove~\eqref{eq:volumecontrol}. Fix some $B= B_i^j$; we can assume  that $B$ intersects $E$ in at least a point $p$, which implies in particular that $B\subset \cBall(E,1/j)$. Without loss of generality, suppose that $p=0$ and $(uf)(p)= 0$; we know that for every  $q\in B\cap \Sigma$
\[
      \dist(q, T^H_0\Sigma) \le \epsilon \Vert q\Vert_\G
      	\qquad\text{and}\qquad 
		\rho_\KK ( u(q) , \DH (uf)_0(q)  ) \le M\epsilon \Vert q\Vert_\G.
\]
Observing that $\DH (uf)_0$ has Lipschitz constant at most $M$, we get
      \begin{align*}
        \dist(u(q), \DH (uf)_0(T_0^H \Sigma)) 
        &\le \rho_\KK ( u(q) , \DH (uf)_0(q)  ) +  M\dist(q , T_0^H \Sigma) \\
&\le 2M \epsilon  \Vert q\Vert_\G.
      \end{align*}
Denoting by $\LL'$ the homogeneous subgroup $\DH (uf)_0(T_0^H \Sigma)$, which is strictly contained in $\LL$, and using the fact that $u(B\cap\Sigma)\subset\LL$, we conclude that 
\[
u(B\cap\Sigma) \subset \cBall_\LL(\LL',2M\epsilon\diam B) \cap \cBall_\LL(0,M\diam B),
\]
where we also used the fact that  the Lipschitz constant of $u|_{B\cap\Sigma}=(uf)|_{B\cap\Sigma}$ is at most $M$. 
By homogeneity one has
\begin{align*}
  \psi^{\ell} (u(B\cap\Sigma))
  &\le (\diam B)^\ell\:\psi^\ell(\cBall_\LL(\LL', 2M\epsilon)\cap\cBall_\LL(0,M))  \\
  &\le M^\ell (\diam B)^\ell\:\psi^\ell(\cBall_\LL(\LL', 2\epsilon)\cap\cBall_\LL(0,1)).
\end{align*}
The claim~\eqref{eq:volumecontrol} follows on  letting
\[
       C(\epsilon,\LL) 
       \defeq \sup_\PP\  \psi^\ell(\cBall_\LL(\PP, 2\epsilon)\cap\cBall_\LL(0,1))        
\]   
where the supremum is taken among proper homogeneous subgroups of $\LL$. 
The fact that $\lim_{\epsilon\to 0^+}C(\epsilon,\LL)=0$ can be easily checked in linear coordinates on the vector space $\LL$,  by comparing $\rho_\LL$ with the Euclidean distance and noting that $\psi^\ell$ is a multiple of the Lebesgue measure.

Combining~\eqref{eq:volumecontrol},~\eqref{eq03261658} and~\eqref{eq:quasiHausdorff}, we  obtain
\begin{align*}
  \int_{\LL} \psi^{Q-m-\ell}_{1/j}(E\cap u^{-1}(s)) \dd \psi^{\ell}(s) \le M^\ell C(\epsilon, \LL) (\psi^{Q-m}(E)+ 1/j)
\end{align*}
and, letting $j\to\infty$, we deduce by Fatou's Lemma that
\[
       \mu_{\Sigma,u} (E) \le M^\ell C(\epsilon,\LL)\psi^{Q-m}(E) .
\]
The proof is accomplished by letting $\epsilon\to 0^+$. 
\end{proof}

Lemma~\ref{lem07241228}, combined with Propositions~\ref{prop:doubling} and~\ref{prop:RadonNykodym}, provides the following consequence.
Recall that $\psi^d$ is Borel regular and that the restriction of a Borel regular measure to a Borel set is Borel regular again.

\begin{Corollary}\label{cor:baddensityzero}
	Under the assumptions and notation of Theorem~\ref{thm:coarea2}, the equality 
	$\Theta(p) = 0$
  	holds for $\psi^{Q-m}$-a.e.~$p\in\Sigma$ such that  $\DH (uf)\vert_{T^H_p \Sigma}$ is not sur\-jec\-ti\-ve on $\LL$. In particular
	\[
	\Theta(p) = \calC(T^H_p\Sigma,\DH (uf)(p)) =0
	\] 
	at all such points $p$.
\end{Corollary}

\subsection{Proof of the coarea formula}
In this section we prove the main coarea formulae of the paper. We start by Theorems~\ref{thm:coarea} and~\ref{thm:coarea2}.

\begin{proof}[Proof of Theorems~\ref{thm:coarea} and~\ref{thm:coarea2}]
Notice that Theorem~\ref{thm:coarea2} implies Theorem~\ref{thm:coarea}.
Sta\-tements $(i)$ and $(ii)$ and the first part of $(iii)$ of Theorem~\ref{thm:coarea2} follow from Proposition~\ref{prop07171752}. The remaining claim~\eqref{eq:coarealocal} follows from Proposition~\ref{prop07251509} and Corollary~\ref{cor:baddensityzero}.
\end{proof}

A direct consequence is Corollary~\ref{cor:coareaPRODUCT}, where we assume that $\KK=\LL\times\MM$ is a direct product:

\begin{proof}[Proof of Corollary~\ref{cor:coareaPRODUCT}]
It is enough to prove the statement in case $R$ is a $C^1_H$ submanifold; actually, we can also assume that there exists $f\in C^1_H(\Omega;\MM)$ such that  $R=\Sigma:=\{p\in\Omega:f(p)=0\}$ and all points in $\Omega$ are split-regular for $f$. Since $\KK=\LL\times\MM$ is a direct product, we have  $uf\in C^1_H(\Omega;\KK)$ and $\DH (uf)_p(g)=\DH  u_p(g)\DH f_p(g)$ for every $g\in\G$. Moreover, since $T^H_p\Sigma=\ker \DH  f_p$,  the equality $\DH (uf)_p|_{T^H_p\Sigma}=\DH  u_p|_{T^H_p\Sigma}$ holds. In particular, condition~\eqref{eq:technicalassumptionPRODUCT} now implies~\eqref{eq:technicalassumption}, and the statement directly follows from Theorem~\ref{thm:coarea}.
\end{proof}

\section{Heisenberg groups}

The most notable examples of Carnot groups are provided by Heisenberg groups. For an integer $n\ge 1$, the $n$-th {\em Heisenberg group} $\HH^n$ is the stratified Lie group associated with the step 2 algebra $V=V_1\oplus V_2$ defined by
\begin{align*}
& V_1=\textrm{span}\{X_1,\dots,X_n,Y_1,\dots,Y_n\},\qquad V_2=\textrm{span}\{T\},\\
& [X_i,Y_j]=\delta_{ij}T\qquad\text{for every }i,j=1,\dots,n.
\end{align*}
We will identify $\HH^n\equiv\R^{2n+1}$ by the {\em exponential coordinates}:
\begin{equation}\nonumber
\R^n\times\R^n\times\R \ni(x,y,t)\longleftrightarrow \exp(x_1 X_1+\dots+y_nY_n+tT)\in\HH^n ,
\end{equation}
according to which the group operation is
\[
(x,y,t)(x',y',t')=(x+x',y+y',t+\tfrac12\textstyle\sum_{j=1}^n(x_jy_j'-x_j'y_j)).
\]
We say that a homogeneous distance $\rho$ on $\HH^n$ is {\em rotationally invariant}\footnote{The terminology ``rotationally invariant'' might be misleading in $\HH^n$ for $n>1$, as not all rotations around the $T$ axis are isometries}  (\cite{SNGRigot}) if
\begin{equation}\label{eq:seba1}
\rho(0,(x,y,t))=\rho(0,(x',y',t))\qquad\text{whenever }|(x,y)|=|(x',y')|,
\end{equation}
where $|\cdot|$ is the Euclidean norm in $\R^{2n}$.
Observe that $\rho$ is rotationally invariant if and only if it is {\it multiradial} according to \cite[Definition 2.21]{CorniMagnani}, i.e., if $\rho(0,(x,y,t))=f(|(x,y)|,|t|)$ for a suitable $f$. 

If $\HH^n=\W\cdot\V$ is a splitting of the $n$-th Heisenberg group $\HH^n$ with $\W$ normal, then necessarily $\V$ is an Abelian {\em horizontal subgroup}, i.e., $\V\subset V_1$, while $\W$ is {\em vertical}, i.e., $V_2\subset \W$. See~\cite[Remark 3.12]{FSSCAdvMath}.
Moreover, if $1\le k\le n$, then the following conditions are equivalent:
\begin{enumerate}[label=(\roman*)]
\item $\PP\subset\HH^n$ is a vertical subgroup  with topological dimension $2n+1-k$;
\item $\PP=P\times V_2$ for some $(2n-k)$-dimensional subspace $P\subset V_1$;
\item $\PP\in\scr T_{\HH^n,\R^k}$.
\end{enumerate}
Proving the equivalence of the statements above is a simple task when one takes into account that every vertical subgroup of codimension at most $n$ possesses a complementary horizontal subgroup, see e.g.~\cite[Lemma~3.26]{FSSCAdvMath}.

\subsection{Area formula in Heisenberg groups}
We provide an explicit representation for the spherical measure on  vertical subgroups of $\HH^n$:

\begin{Proposition}\label{prop:Sverticalsubgroups}
Assume that $\HH^n$ is endowed with a rotationally invariant homogeneous distance and let  $1\le k\le n$. Then, there exists a constant $c(n,k)$ such that for every vertical subgroup $\PP\in\scr T_{\HH^n,\R^k}$
\[
c(n,k)\cal S^{2n+2-k}\hel \PP = \cal H^{2n+1-k}_E\hel\PP,
\]
where $\cal H^{2n+1-k}_E$ denotes the Euclidean Hausdorff measure on $\R^{2n+1}\equiv\HH^n$.
\end{Proposition}
\begin{proof}
Let $\PP\in\scr T_{\HH^n,\R^k}$ be a fixed vertical subgroup; by~\cite[Lemma~3.26]{FSSCAdvMath}  there exists a complementary Abelian horizontal subgroup $\V = V\times\{0\}$, for a proper $k$-dimensional subspace $V\subset V_1$.
Let $W$ be a $(2n-k)$-dimensional complementary subspace of $V$ in $V_1$ and set $\W:=W\times V_2$, which is a vertical subgroup that is complementary to $\V$.
Let $P\subset V_1$ such that $\PP=P\times V_2$.

Let $f:W\to V$ be such that $P=\{w+f(w):w\in W\}$ and let $\phi:\W\to V$ such that $\PP=\{w(\phi(w),0):w\in\W\}$.
Now, notice that if $z\in W$ and $t\in\R$, then 
\[
(z,t)(\phi(z,t),0)=(z+\phi(z,t),t+\frac12\omega(z,\phi(z,t)) ,
\]
where $\omega$ is the standard symplectic form on $\R^{2n}$.
Since $z+\phi(z,t)\in P$ and $(z+V)\cap P = \{z+f(z)\}$, then we have $\phi(z,t)=f(z)$.

The area formula of \cite[Theorem~1.2]{CorniMagnani}, together with \cite[Theorem~2.12 and Proposition 2.13]{CorniMagnani} from the same paper, provide a constant $c(n,k)>0$ such that
\begin{equation}\label{eq:uno11}
c(n,k)\cal S^{2n+2-k}\hel \PP = \Phi_\#(J^\phi\phi\:\cal H^{2n+1-k}_E\hel\W),
\end{equation}
where $J^\phi\phi$ is the intrinsic Jacobian of $\phi$ as in \cite[Definition~2.14]{CorniMagnani} and $\Phi$ is the intrinsic graph map. On the other side, the Euclidean area formula gives 
\begin{equation}\label{eq:due22}
\cal H^{2n+1-k}_E\hel \PP = F_\#(JF\:\cal H^{2n+1-k}_E\hel\W),
\end{equation}
where $F:\W\to\PP$ is defined by $F(x,y,t):=(f(x,y),t)$ for every $(x,y)\in W$ and 
$JF$ is the Euclidean area factor. As a matter of fact, using the equality $f=\phi$, one has $J^\phi\phi=JF$ and the statement immediately follows from~\eqref{eq:uno11} and~\eqref{eq:due22}.
\end{proof}

\begin{Remark}
Proposition~\ref{prop:Sverticalsubgroups} holds, with no changes in the proof, in the more general case $\HH^n$ is endowed with a homogeneous distance that is {\it $(2n+1-k)$-vertically symmetric} according to~\cite[Definition~2.19]{CorniMagnani}. 
\end{Remark}

\begin{Remark}\label{rem:Hrotatinvar}
	When $\HH^n$ is endowed with a rotationally invariant distance $\rho$, then for every pair $(\PP,\PP')$ of one-codimensional homogeneous subgroups of $\HH^n$, 
	there exist an isometry $(\HH^n,\rho)\to(\HH^n,\rho)$ that maps $\PP$ to $\PP'$.
	The proof is left  to the reader.
\end{Remark}

The following proposition completes the proof of Corollary~\ref{cor:SvsH}.
\begin{Proposition}\label{prop:Hncodim1}
	If $\HH^n$ is endowed with a rotationally invariant homogeneous distance and $\G'=\R$, then the function $\textfrak a$ in Corollary \ref{cor:SvsH} is constant, i.e., there exists $C\in[1,2^{2n+1}]$ such that 
\[
\text{$\cal S^{2n+1} \hel R= C\cal H^{2n+1}\hel R\qquad \forall\:(\HH^n,\R)$-rectifiable set $R\subset\HH^n$.}
\]
\end{Proposition}
\begin{proof}
	When $\G=\HH^n$ and $\G'=\R$, then the function $\textfrak a$ defined in~\eqref{eq03231608} is constant by Remark~\ref{rem:Hrotatinvar}.
\end{proof}

Similarly, Corollary~\ref{cor:densityexistence} can be improved when $\G$ is the Heisenberg group endowed with a rotationally invariant distance.

\begin{Corollary}\label{cor:densityconstantinHeis}
Assume $\G$ is the Heisenberg group $\HH^n$ endowed with a rotationally invariant distance and $\G'=\R^m$ for some $1\le m\le n$; if $\psi^{2n+2-m}$ is the spherical Hausdorff measure, then the function $\textfrak d$ in Corollary~\ref{cor:densityexistence} is constant.

If $m=1$ and $\psi^{2n+2-m}$ is the  Hausdorff measure, then the function $\textfrak d$ in Corollary~\ref{cor:densityexistence} is constant.
\end{Corollary}
\begin{proof}
Concerning the first part of the statement, let $\W\in\scr T_{\HH^n,\R^m}$ be fixed; by Proposition~\ref{prop:Sverticalsubgroups} we have
\begin{align*}
\textfrak d(\W)&=\lim_{r\to 0^+} \frac{\cal S^{2n+2-m}(\W\cap\Ball(0,r))}{r^{2n+2-m}}\\
& = \cal S^{2n+2-m}(\W\cap\Ball(0,1)) =c(n,m) \cal H_E^{2n+1-m}(\W\cap\Ball(0,1))
\end{align*}
and the latter quantity does not depend on $\W$ by  rotational invariance of the distance. 
The second part of the statement is an immediate consequence of Remark~\ref{rem:Hrotatinvar}.
\end{proof}

\subsection{Coarea formula in Heisenberg groups}
When one considers spherical measures in the Heisenberg group endowed with a rotationally invariant distance,
then  the coarea factor coincides up to a multiplicative constant with the quantity
\begin{equation*}
\begin{split}
& J^Ru(p):=\left( \det (L\circ L^T) \right)^{1/2},\qquad L:=\DH u_p|_{T^H_pR}.
\end{split}
\end{equation*}
We prove this fact.

\begin{Proposition}\label{prop:coareafactorH}
Consider the Heisenberg group $\HH^n$ endowed with a rotationally invariant distance. Let $\PP\in\scr T_{\HH^n,\Rm}$ be a vertical subgroup of topological dimension $2n+1-m$ and let $L:\PP\to\R^\ell$ be a homogeneous morphism; assume $1\le m+\ell\le n$. Then
\[
\calC(\PP,L)=\frac{c(n,m+\ell)}{c(n,m)}\left( \det (L\circ L^T) \right)^{1/2},
\]
where the positive constants $c(n,m)$ and $c(n,m+\ell)$ are those provided by Proposition~\ref{prop:Sverticalsubgroups}.
\end{Proposition}
\begin{proof}
 If $L$ is not onto $\R^\ell$, then the statement is true. 
 We assume that $L$ is surjective. 
 By Proposition~\ref{prop:Sverticalsubgroups}
\begin{align*}
\mu_{\PP,L} &=\int_{\R^\ell}\cal S^{2n+2-m-\ell}\hel L^{-1}(s)\dd\cal L^\ell(s)\\
& = c(n,m+\ell)\int_{\R^\ell} \cal H^{2n+1-m-\ell}_E\hel L^{-1}(s)\dd\cal L^\ell(s)\\
& = c(n,m+\ell) (\det  L\circ L^T)^{1/2}\cal H^{2n+1-m}_E\hel \PP,
\end{align*}
where we used the Euclidean coarea formula. A second application of Proposition~\ref{prop:Sverticalsubgroups} gives 
\[
\mu_{\PP,L} = \frac{c(n,m+\ell)}{c(n,m)} (\det  L\circ L^T)^{1/2}\cal S^{2n+2-m}\hel \PP
\]
and this is enough to conclude.
\end{proof}

We now have all the tools needed in order to prove our coarea formula in Heisenberg groups.

\begin{proof}[Proof of Theorem~\ref{thm:coareaHeisenberg}]
The first part of the statement is an immediate consequence of Corollary~\ref{cor:coareaPRODUCT} and the fact that, if $\DH  u_p|_{T^H_pR}$ is surjective on $\R^\ell$, then $T^H_pR\cap\ker\DH  u_p$ is a vertical subgroup of dimension $2n+1-m-\ell\ge n+1$, and  by~\cite[Lemma~3.26]{FSSCAdvMath} it admits a complementary (horizontal) subgroup.

The second part of the statement is now a consequence of Proposition~\ref{prop:coareafactorH}; clearly, one has $\textfrak c=c(n,m+\ell)/c(n,m)$ according to the constants introduced in Proposition~\ref{prop:Sverticalsubgroups}.
\end{proof}

\bibliographystyle{abbrv}
\bibliography{RefsCarnot}

\end{document}